\documentclass[a4paper]{article}
\usepackage[latin1]{inputenc}
\usepackage[T1]{fontenc}
\usepackage[francais,english]{babel}
\usepackage{amsmath}
\usepackage{amssymb}
\usepackage{theorem}
\newtheorem{thm}{Theorem}
\newtheorem{Prop}{Proposition}
\newtheorem{Def}{Definition}
\newtheorem{Lem}{Lemma}
\newtheorem{Cor}{Corollary}
\newtheorem{rk}{Remark}

\author{
Karine \textsc{Beauchard}\footnote{
CMLA, ENS Cachan, CNRS, UniverSud, 61, avenue du Président Wilson, F-94230 Cachan, FRANCE.
email: Karine.Beauchard@cmla.ens-cachan.fr}
\thanks{The author was partially supported by the ``Agence Nationale de la Recherche'' (ANR),
Projet Blanc C-QUID number BLAN-3-139579}, 
Camille \textsc{Laurent}\footnote{
Laboratoire Mathématiques, UMR 8628 du CNRS, Bâtiment 425, Faculté des Sciences d'Orsay, 
Université Paris-Sud 11, F-91405 Orsay Cedex, France.
email: camille.laurent@math.u-psud.fr }
}

\title{Local controllability of 1D linear and nonlinear Schr\"odinger equations 
with bilinear control}
\date{}

\begin{document}
\maketitle

\begin{abstract}
We consider a linear Schr\"odinger equation, on a bounded interval, 
with bilinear control, that represents a quantum particle in 
an electric field (the control). We prove the controllability 
of this system, in any positive time, locally around the ground state.

Similar results were proved for particular models  
\cite{KB-JMPA, SchroLgVar, KB-JMC},
in non optimal spaces, in long time 
and the proof relied on the Nash-Moser implicit function theorem 
in order to deal with an a priori loss of regularity.

In this article, the model is more general,
the spaces are optimal, there is no restriction on the time 
and the proof relies on the classical inverse mapping theorem. 
A hidden regularizing effect is emphasized,
showing there is actually no loss of regularity.

Then, the same strategy is applied to nonlinear Schr\"odinger equations
and nonlinear wave equations, showing that the method works for a wide range
of bilinear control systems.

\begin{center}
\textbf{Résumé}
\end{center}

On considère une équation de Schr\"odinger linéaire, sur un intervalle borné,
avec contrôle bilinéaire, représentant une particule quantique dans un champ
électrique (le contrôle). On démontre la contrôlabilité locale de ce système,
en tout temps positif, localement au voisinage de l'état fondamental.

Des résultats similaires ont déjà été établis \cite{KB-JMPA, SchroLgVar, KB-JMC},
mais dans des espaces non optimaux, en temps long et leur preuve reposait
sur le théorème de Nash-Moser, pour gérer une apparente perte de régularité.

Dans cet article, le modèle étudié est plus général, les espaces 
sont optimaux, il n'y a pas de restriction sur le temps et la preuve repose sur
le théorème d'inversion locale classique. Un effet régularisant est
exhibé, montrant qu'il n'y a finalement pas de perte de régularité.

La même stratégie est ensuite utilisée sur des équations de Schr\"odinger 
nonlinéaires et des équations des ondes nonlinéaires, montrant qu'elle
s'applique de façon assez générale aux systèmes de contrôle bilinéaires.
\end{abstract}


\section{Introduction}

\subsection{Main result}

Following \cite{rouchon-oct2002}, we consider a quantum particle, 
in a 1D infinite square potential well, subjected to an electric field. 
It is represented by the following Schr\"odinger equation
\begin{equation} \label{Schro_eq}
\left\lbrace\begin{array}{l}
i \frac{\partial \psi}{\partial t}(t,x) = 
- \frac{\partial^2 \psi}{\partial x^2}(t,x)
- u(t) \mu(x) \psi (t,x) , x \in (0,1), t \in (0,T), \\
\psi(t,0)=\psi(t,1)=0,
\end{array}\right.
\end{equation}
where $\psi$ is the wave function of the particle,
$u$ is the amplitude of the electric field and
$\mu \in H^{3}((0,1),\mathbb{R})$ is the dipolar moment.
The system (\ref{Schro_eq}) is a bilinear control system, in which
\begin{itemize}
\item the state is $\psi$, with $\| \psi(t) \|_{L^2(0,1)} = 1$, $\forall t \in (0,T)$,
\item the control is the real valued function $u:[0,T] \rightarrow \mathbb{R}$.
\end{itemize}

Let us introduce some notations. The operator $A$ is defined by
\begin{equation} \label{def:A}
\begin{array}{cc}
D(A):=H^2 \cap H^1_0((0,1),\mathbb{C}), 
&
A \varphi := - \frac{d^2 \varphi}{dx^2}.
\end{array}
\end{equation}
Its eigenvalues and eigenvectors are
\begin{equation} \label{vap_vep}
\lambda_k:=(k\pi)^2, 
\varphi_k(x):=\sqrt{2} \sin (k \pi x), \forall k \in \mathbb{N}^*.
\end{equation}
The family $(\varphi_k)_{k \in \mathbb{N}^*}$ is an orthonormal basis
of $L^2((0,1),\mathbb{C})$ and
$$\psi_k(t,x) := \varphi_k(x) e^{-i\lambda_k t}, \forall k \in \mathbb{N}^*$$
is a solution of (\ref{Schro_eq}) with $u \equiv 0$ called eigenstate,
or ground state, when $k =1$. We define the spaces
\begin{equation} \label{def:Hs}
H^s_{(0)}((0,1),\mathbb{C}):=D(A^{s/2}), \forall s>0
\end{equation}
equipped with the norm
$$\| \varphi \|_{H^s_{(0)}} := \left( \sum\limits_{k=1}^{\infty}
| k^s \langle \varphi , \varphi_k \rangle |^2 \right)^{1/2}.$$
We denote by $\langle . , \rangle$ the $L^2((0,1),\mathbb{C})$
scalar product
$$\langle f , g \rangle = \int_0^1 f(x) \overline{g(x)} dx$$
and by $\mathcal{S}$ the unit $L^2((0,1),\mathbb{C})$-sphere.
The first goal of this article is the proof of the following result.

\begin{thm} \label{Main-Thm}
Let $T>0$ and $\mu \in H^{3}((0,1),\mathbb{R})$ be such that 
\begin{equation} \label{hyp_mu}
\exists c>0 \text{ such that } 
\frac{c}{k^3} \leqslant |\langle \mu \varphi_1 , \varphi_k \rangle |, 
\forall k \in \mathbb{N}^*.
\end{equation}
There exists $\delta>0$ and a $C^1$ map 
$$\begin{array}{cccl}
\Gamma: & \mathcal{V}_T & \rightarrow & L^2((0,T),\mathbb{R})
\end{array}$$
where
$$\mathcal{V}_T := \{ \psi_f \in \mathcal{S} \cap H^3_{(0)}((0,1),\mathbb{C}) ;
\| \psi_f - \psi_1(T) \|_{H^3} < \delta  \},$$
such that, $\Gamma( \psi_1(T) ) =0$ and for every $\psi_f \in  \mathcal{V}_T$,
the solution of (\ref{Schro_eq}) with initial condition 
\begin{equation} \label{IC}
\psi(0)=\varphi_1 
\end{equation}
and control $u=\Gamma(\psi_f)$ satisfies $\psi(T)=\psi_f$.
\end{thm}

\begin{rk}
Thanks to the time reversibility of the system,
Theorem \ref{Main-Thm} ensures the local controllability of the system (\ref{Schro_eq})
around the ground state: for every $T>0$, there exists $\delta>0$ such that,
for every $\psi_0, \psi_f \in \mathcal{S} \cap H^3_{(0)}((0,1),\mathbb{C})$
with $\|\psi_0-\psi_1(0)\|_{H^3}+\| \psi_f-\psi_1(T)\|_{H^3} < \delta$,
there exists a control $u \in L^2(0,T)$ such that the solution of 
(\ref{Schro_eq}) with initial condition $\psi(0)=\psi_0$ satisfies $\psi(T)=\psi_f$.
\end{rk}

\begin{rk}
The assumption (\ref{hyp_mu}) holds, for example, with $\mu(x)=x^{2}$, because
\begin{equation} \label{coeff-explicit-x2}
\langle x^{2} \varphi_{1},\varphi_{k} \rangle = 
\int_{0}^{1} 2 x^{2} \sin (k\pi x) \sin(\pi x) dx =
\left\lbrace \begin{array}{l}
\frac{(-1)^{k+1} 8k}{\pi^{2} (k^{2}-1)^{2}} \text{ if } k \geqslant 2,
\\
\frac{-3+2\pi^{2}}{6\pi^{2}} \text{ if } k=1.
\end{array}\right.
\end{equation}
But it does not hold when $\langle \mu \varphi_{1} , \varphi_{k} \rangle = 0$,
for some $k \in \mathbb{N}^{*}$,
or when $\mu$ has a symmetry with respect to $x=1/2$.
However, the assumption (\ref{hyp_mu}) holds
generically with respect to $\mu \in H^{3}((0,1),\mathbb{R})$ because
\begin{equation} \label{hyp_mu-generic}
\langle \mu \varphi_1 , \varphi_k \rangle 
=
\frac{4[(-1)^{k+1} \mu'(1) - \mu'(0)]}{k^3 \pi^2}
- \frac{\sqrt{2}}{(k \pi)^3} \int_0^1 ( \mu \varphi_1)'''(x) \cos(k \pi x ) dx,
\forall k \in \mathbb{N}^*.
\end{equation}
(see Appendix \ref{appendix:hyp_mu} for a proof). Thus, Theorem 1 is very general.
\end{rk}

\subsection{A simpler proof}

The local controllability of 1D Schr\"odinger equations, with bilinear control,
has already been investigated in 
\cite{KB-JMPA, SchroLgVar, KB-JMC},
(see also \cite{KB-poutres} for a similar result on a 1D beam equation).
In these articles, three different models are studied. The local controllability
of the nonlinear system is proved thanks to the linearization principle:
\begin{itemize}
\item first, we prove the controllability of a linearized system,
\item then, we prove the local controllability of the nonlinear system, 
by applying an inverse mapping theorem.
\end{itemize}
This strategy is coupled with the return method and quasi-static deformations
in \cite{KB-JMPA, KB-JMC}  and with power series expansions in 
\cite{SchroLgVar, KB-JMC} (see \cite{jmc-tank, JMC-book} by Coron
for a presentation of these technics). 
In these articles, the most difficult part of the proof is
the application of the inverse mapping theorem. Indeed, because of an a priori loss
of regularity, we were lead to apply the Nash-Moser implicit function theorem
(see, for instance \cite{Alinhac-Gerard} by Alinhac, Gérard and \cite{horm} by H\"ormander), 
instead of the classical inverse mapping theorem. The Nash-Moser
theorem requires, in particular, the controllability of an infinite number of
linearized systems, and tame estimates on the corresponding controls. These two points
are difficult to prove and lead to long technical developments in 
\cite{KB-JMPA, SchroLgVar, KB-JMC}.

In this article, we propose a simpler proof, that uses only the classical 
inverse mapping theorem (needing the controllability of only one linearized system), 
because we emphasize a hidden regularizing effect (see Proposition \ref{WP-CYpb}). 

Therefore, the controllability result of Theorem \ref{Main-Thm} enters 
the classical framework of local controllability results for nonlinear systems,
proved with fixed point arguments 
(see, for instance, 
\cite{RosierKdV} by Rosier, 
\cite{Crepeau_Cerpa} by Cerpa and Crépeau,
\cite{Russel_Zhang428} by Russell and Zhang, 
\cite{Zhang_513} by Zhang,
\cite{zuazua-wave} by Zuazua; this list is not exhaustive).

\subsection{Additionnal results}

The proof we developed for Theorem \ref{Main-Thm}
is quite robust, thus we could apply it to other situations:
other linear PDEs and also nonlinear PDEs,
that are presented in the next subsections.
This shows that the strategy proposed in this article 
works for a wide range of bilinear systems.

\subsubsection{Generalization to higher regularities}

The first situation is the analogue result of Theorem \ref{Main-Thm},
but with higher regularities: we prove 
the local exact controllability of (\ref{Schro_eq}) in smoother spaces 
and with smoother controls. Namely, we prove the following result.

\begin{thm} \label{Main-ThmH5H1}
Let $T>0$ and $\mu \in H^{5}((0,1),\mathbb{R})$ be such that 
(\ref{hyp_mu}) holds. There exists $\delta>0$ and a $C^1$ map 
$$\begin{array}{cccl}
\Gamma: & \mathcal{V}_T & \rightarrow & H^1_0((0,T),\mathbb{R})
\\
        &   \psi_f      & \mapsto     & \Gamma(\psi_f)
\end{array}$$
where
$$\mathcal{V}_T := \{ \psi_f \in \mathcal{S} \cap H^5_{(0)}((0,1),\mathbb{C}) ;
\| \psi_f - \psi_1(T) \|_{H^5} < \delta  \},$$
such that, $\Gamma( \psi_1(T) ) =0$ and for every $\psi_f \in  \mathcal{V}_T$,
the solution of (\ref{Schro_eq}), (\ref{IC}) 
with control $u=\Gamma(\psi_f)$ satisfies $\psi(T)=\psi_f$.
\end{thm}

Of course, the strategy may be used to go further and prove the local
exact controllability of (\ref{Schro_eq}) around the ground state
\begin{itemize}
\item in $H^7_{(0)}(0,1)$ with controls in $H^2_0((0,T),\mathbb{R})$, 
\item in $H^9_{(0)}(0,1)$ with controls in $H^3_0((0,T),\mathbb{R})$, etc.
\end{itemize}

\subsubsection{On the 3D ball with radial data}

The second situation is the analogue result of Theorem \ref{Main-Thm},
but for the Schr\"odinger equation posed on the three dimensional unit ball $B^3$
for radial data. In polar coordinates, the Laplacian for radial data can be written
$$\Delta u(r)=\partial_r^2 u(r)+ \frac{2}{r}\partial_r u(r).$$
In particular, we have $\Delta \left(\frac{g(r)}{r}\right)=\frac{\partial_r^2 u(r)}{r}$. 
The eigenfunctions of the Dirichlet operator $A=-\Delta$ with domain 
$D(A):=H^2_{radial} \cap H^1_0(B^3)$ are 
$\varphi_k=\frac{sin(k\pi r)}{r\sqrt{2\pi}}$ with eigenvalues 
$\lambda_k=(k\pi)^2$. Thus, we study the Schr\"odinger equation
\begin{equation} \label{Schro_eq_rad}
\left\lbrace \begin{array}{l}
i \frac{\partial \psi}{\partial t}(t,r) =
- \Delta \psi(t,r) - u(t) \mu(r) \psi(t,r), r \in (0,1), \\
\psi(t,1)=0.
\end{array} \right.
\end{equation}

The theorem we obtain is very similar to Theorem \ref{Main-Thm}.

\begin{thm} \label{radial-Thm}
Let $T>0$ and $\mu \in H^{3}(B^3,\mathbb{R})$ radial be such that 
\begin{equation} \label{hyp_mu_rad}
\exists c>0 \text{ such that } 
\frac{c}{k^3} \leqslant |\langle \mu \varphi_1 , \varphi_k \rangle |, 
\forall k \in \mathbb{N}^*.
\end{equation}
There exists $\delta>0$ and a $C^1$ map 
$$\begin{array}{cccl}
\Gamma: & \mathcal{V}_T & \rightarrow & L^2((0,T),\mathbb{R})
\end{array}$$
where
$$\mathcal{V}_T := \{ \psi_f \in \mathcal{S} \cap H^3_{(0),rad}(B^3,\mathbb{C}) ;
\| \psi_f - \psi_1(T) \|_{H^3} < \delta  \},$$
such that, $\Gamma( \psi_1(T) ) =0$ and for every $\psi_f \in  \mathcal{V}_T$,
the solution of (\ref{Schro_eq_rad}) with initial condition 
\begin{equation} \label{IC_rad}
\psi(0)=\varphi_1 
\end{equation}
and control $u=\Gamma(\psi_f)$ satisfies $\psi(T)=\psi_f$.
\end{thm}

The analysis is very close to the 1D case since for this particular data, 
the Laplacian behaves as in dimension 1. We refer to Appendix A for the proof of 
the genericity of the assumption (\ref{hyp_mu_rad}). Note that this simpler situation 
has also been used by Anton for proving global 
existence for the nonlinear Schr\"odinger equation \cite{AntonNLSbouleradial}.

\subsubsection{Nonlinear Schr\"odinger equations}

The third situation concerns nonlinear Schr\"odinger equations.
More precisely we study the following nonlinear Schr\"odinger equation
with  Neumann boundary conditions
\begin{equation} \label{NLS}
\left\lbrace \begin{array}{l}
i \frac{\partial \psi}{\partial t}(t,x) = - \frac{\partial^2 \psi}{\partial x^2}(t,x)
+ |\psi|^2 \psi(t,x) - u(t) \mu(x) \psi(t,x), x \in (0,1), t \in (0,T),\\
\frac{\partial \psi}{\partial x}(t,0)=\frac{\partial \psi}{\partial x}(t,1)=0.
\end{array}\right.
\end{equation}
It is a nonlinear control system where
\begin{itemize}
\item the state is $\psi$, with $\|\psi(t)\|_{L^2(0,1)}=1, \forall t \in [0,T]$,
\item the control is the real valued function $u:[0,T]\rightarrow \mathbb{R}$.
\end{itemize}
We study its local controllability around the reference trajectory
$$(\psi_{ref}(t,x):=e^{-it},u_{ref}(t)=0).$$
More precisely, we prove the following result.

\begin{thm}\label{thm_controlNLS}
Let $T>0$ and $\mu \in H^2(0,1)$ be such that
\begin{equation} \label{hyp_mu-NLS}
\exists c>0  \text{  such that  }  
\Big| \int_0^1 \mu(x) \cos( k \pi x) dx \Big|
\geqslant \frac{c}{\max\{1,k\}^2},
\forall k \in \mathbb{N}.
\end{equation}
There exists $\eta>0$ and a $C^1$-map
$$\Gamma : \mathcal{V}_T \rightarrow L^2((0,T),\mathbb{R})$$
where
$$\mathcal{V}_T := \{ \psi_f \in \mathcal{S} \cap H^2 (0,1) ;
\psi_f'(0)=\psi_f'(1)=0 \text{ and }
\| \psi_f - e^{-iT} \|_{H^2} < \eta \}$$
such that, for every $\psi_f \in \mathcal{V}_T$, 
the solution of (\ref{NLS}) with initial condition
\begin{equation} \label{IC-NLS}
\psi(0,x)=1, \forall x \in (0,1)
\end{equation}
and control $u:=\Gamma(\psi_f)$ is defined on $[0,T]$ and satisfies $\psi(T)=\psi_f$.
\end{thm}

\begin{rk}
The assumption (\ref{hyp_mu-NLS}) holds generically in $H^2(0,1)$.
Indeed, integrations by part give
$$\int_0^1 \mu(x) \cos( k \pi x) dx
=\frac{1}{(k\pi)^2} \left(
(-1)^{k+1} \mu'(1) + \mu'(0) + \int_0^1 \mu''(x) \cos(k\pi x) dx 
\right), \forall k \in \mathbb{N}^*.$$
\end{rk}

Other versions of this result, with higher regularities may be proved:
the system is exactly controllable, locally around the reference trajectory
\begin{itemize}
\item in $H^4(0,1)$ with controls in $H^1_0(0,T)$,
\item in $H^6(0,1)$ with controls in $H^2_0(0,T)$, etc.
\end{itemize}
Focusing nonlinearities may also be considered.

\subsubsection{Nonlinear wave equations}

The third situation concerns nonlinear wave equations.
More precisely we study the following wave equation
with  Neumann boundary conditions
\begin{equation} \label{wave}
\left\lbrace \begin{array}{l}
w_{tt}=w_{xx}+f(w,w_t)+u(t) \mu(x) (w+ w_t), x \in (0,1), t \in (0,T),\\
w_x(t,0)=w_x(t,1)=0,
\end{array}\right.
\end{equation}
where $f$ is an appropriate nonlinearity, that satisfies, in particular, $f(1,0)=0$. 
It is a nonlinear control system where
\begin{itemize}
\item the state is $(w,w_t)$,
\item the control is the real valued function $u:[0,T] \rightarrow \mathbb{R}$.
\end{itemize}
We study its exact controllability, locally around the reference trajectory
$$( w_{ref}(t,x)=1 , u_{ref}(t)=0 ).$$
More precisely, we prove the following result.

\begin{thm} \label{thm:ondes}
Let $T>2$, $\mu \in H^2((0,1),\mathbb{R})$ be such that (\ref{hyp_mu-NLS}) holds
and $f \in C^3(\mathbb{R}^2,\mathbb{R})$ be such that
$f(1,0)=0$ and $\nabla f (1,0)=0$.
There exists $\eta>0$ and a $C^1$-map
$$\Gamma : \mathcal{V}_T \rightarrow L^2((0,T),\mathbb{R})$$
where
$$\begin{array}{ll}
\mathcal{V}_T := \{ (w_f,\dot{w}_f) \in H^3 \times H^2((0,1),\mathbb{R}) ;
& w_f'(0)=w_f'(1)=\dot{w}_f'(0)=\dot{w}_f'(1)=0  
\\ &
\text{ and } \| w_f - 1 \|_{H^3} + \| \dot{w}_f \|_{H^2} < \eta \}
\end{array}$$
such that $\Gamma(1,0)=0$ and for every $(w_f,\dot{w}_f) \in \mathcal{V}_T$, 
the solution of (\ref{wave}) with initial condition
\begin{equation} \label{IC-wave}
(w,w_t)(0,x)=(1,0), \forall x \in (0,1)
\end{equation}
and control $u:=\Gamma(w_f,\dot{w}_f)$ is defined on $[0,T]$ 
and satisfies $(w,w_t)(T)=(w_f,\dot{w}_f)$.
\end{thm}

Other versions of this result, with higher regularities may be proved:
the system is exactly controllable, locally around the reference trajectory
\begin{itemize}
\item in $H^4 \times H^3(0,1)$ with controls in $H^1_0(0,T)$,
\item in $H^5 \times H^4(0,1)$ with controls in $H^2_0(0,T)$, etc.
\end{itemize}

\subsection{A brief bibliography}

\subsubsection{A previous negative result}
\label{subsub:neg}

First, let us recall an important negative controllability result,
for the equation (\ref{Schro_eq}), proved by Turinici \cite{Defranceschi-LeBris}.
It is a corollary of a more general result due to
Ball, Marsden and Slemrod \cite{ball-marsden-slemrod}.

\begin{Prop} \label{Prop:BMS}
Let $\psi_0 \in \mathcal{S} \cap H^2_{(0)}((0,1),\mathbb{C})$ and 
$U[T;u,\psi_0]$ be the value at time $T$
of the solution of (\ref{Schro_eq}) with initial condition $\psi(0)=\psi_0$.
The set of attainable states from $\psi_0$,
$$\{ U[T;u,\psi_0] ; T>0, u \in L^2((0,T),\mathbb{R}) \}$$
has an empty interior in $\mathcal{S} \cap H^2_{(0)}((0,1),\mathbb{C})$.
Thus (\ref{Schro_eq}) is not controllable in 
$\mathcal{S} \cap H^2_{(0)}((0,1),\mathbb{C})$
with controls in $L^2_{loc}([0,+\infty),\mathbb{R})$.
\end{Prop}

Proposition \ref{Prop:BMS} is a rather weak negative controllability result,
because it does not prevent from positive controllability results,
in different spaces. This had already been emphasized for the particular
cases studied in \cite{KB-JMPA, SchroLgVar, KB-JMC}, 
in which the reachable set is proved to contain $H^{7}_{(0)}$ or
$H^{5+}_{(0)}$.
In this article, we prove that the reachable set (at least locally, with small
controls in $L^2((0,T),\mathbb{R})$), coincides with $\mathcal{S} \cap H^3_{(0)}$,
(which has, indeed, an empty interior in $\mathcal{S} \cap H^2_{(0)}$).
Therefore, sometimes, Ball Marsden and Slemrod's negative result is only due to an 
'unfortunate' choice of functional spaces, that does not allow the controllability. 
It may not be due to a deep non controllability (such as, for example, 
when a subsystem evolves independently of the control).

\subsubsection{Iterated Lie brackets}

Now, let us quote some articles about the controllability of quantum systems.

First, the controllability of \emph{finite dimensional} quantum systems 
(i.e. modelled by an ordinary differential equation) is well understood.
Let us consider the quantum system
\begin{equation} \label{df}
i \frac{dX}{dt}= H_{0}X + u(t) H_{1}X,
\end{equation}
where $X \in \mathbb{C}^{n}$ is the state, 
$H_{0}, H_{1}$ are $n*n$ hermitian matrices, 
and $t \mapsto u(t) \in \mathbb{R}$ is the control.
The controllability of (\ref{df}) is linked to the rank of the Lie algebra 
spanned by $H_{0}$ and $H_{1}$ (see for instance 
\cite{Albertini} by Albertini and D'Alessandro, 
\cite{Altafini} by Altafini, 
\cite{Brockett} by Brockett,
see also \cite{Agrachev-book} by Agrachev and Sachkov,
\cite{JMC-book} by Coron for a more general discussion). 
\\

In \emph{infinite dimension}, there are cases where the iterated Lie brackets 
provide the right intuition. For instance, it  holds for the non controllability of the
harmonic oscillator (see \cite{HarmOsc} by Mirrahimi and Rouchon). 
However, the Lie brackets are often 
less powerful in infinite dimension than in finite dimension. It is precisely
the case of our system. Indeed, let us define the operators
$$\begin{array}{ll}
D(f_{0}):=H^{2} \cap H^{1}_{0}(0,1), 
&
f_{0}(\psi):= -\psi '',
\\
D(f_{1}):= L^{2}(0,1),
&
f_{1}(\psi):= x^{2} \psi,
\end{array}$$
which correspond to $\mu(x)=x^2$.
Let us compute the iterated Lie brackets at the point 
$\varphi_{1}(x)=\sqrt{2} \sin(\pi x)$.
Since $\varphi_{1} \in D(f_{0})$, we can compute
$$\begin{array}{c}
\lbrack f_{0},f_{1} \rbrack (\varphi_{1})= - 4 x \varphi_{1}' - 2 \varphi_{1}, \\
\lbrack f_{1}, \lbrack f_{0},f_{1} \rbrack  \rbrack (\psi)= 
8 x^{2} \varphi_{1} = 8 f_{1}(\varphi_{1}).
\end{array}$$
Notice that $\lbrack f_{0},f_{1} \rbrack (\varphi_{1})$ does not belong to
$D(f_{0})$ because 
$\lbrack f_{0},f_{1} \rbrack (\varphi_{1})(1)= 4 \sqrt{2} \pi \neq 0$.
Thus, in order to give a sense to the Lie bracket
$[f_{0},[f_{0},f_{1}]]$, one needs to extend the definition of $f_{0}$ to
functions that do not vanish at $x=0,1$. A natural choice is
\begin{equation} \label{f0-newdef}
f_{0}(\psi):=-\psi'' + \psi(0) \delta_{0}' - \psi(1) \delta_{1}'
\end{equation}
because, with this choice, we have
$$\langle f_{0}(\psi) , \widetilde{\psi} \rangle = 
\langle \psi , f_{0}(\widetilde{\psi}) \rangle ,
\forall \psi \in D(f_{0}) , 
\forall \widetilde{\psi} \in H^{2}(0,1),$$
in the sense
$$- \int_{0}^{1} \psi''(x) \widetilde{\psi}(x) dx =
- \int_{0}^{1} \psi(x) \widetilde{\psi}''(x) dx 
- \psi'(1) \widetilde{\psi}(1) + \psi'(0) \widetilde{\psi}(0).$$
With the definition (\ref{f0-newdef}), we get
$$[f_{0},[f_{0},f_{1}]](\psi)=-8f_{0}(\psi) + 4 \psi'(1) \delta'_{1}$$
But then, again, $[f_{0},[f_{0},[f_{0},f_{1}]]]$ is not well defined.
Moreover, even if we could give a sense to any iterated Lie bracket, 
because of the presence of Dirac masses, it would not be clear which space 
the Lie algebra should generate in case of local controllability.
Therefore, the way the Lie algebra rank condition could be used directly in
infinite dimension is not clear (see also \cite{JMC-book} for the same discussion
on other examples). This is why we develop completely analytic methods in this article.
\\
\\

Finally, let us quote important articles about the controllability of PDEs,
in which positive results are proved by applying geometric control methods
to the (finite dimensional) Galerkin approximations of the equation.
In \cite{Agrachev9} by Sarychev and Agrachev and 
\cite{Shirikyan446} by Shirikyan, the authors prove exact controllability results
for dissipative equations. In \cite{Chambrion-et-al}, 
by Boscain, Chambrion, Mason and Sigalotti, 
the authors prove the approximate controllability in $L^2$, 
for bilinear Schr\"odinger equations such as (\ref{Schro_eq}).
\\

We also refer to the following works about the controllability of
finite dimensional quantum systems
\cite{Agrachev-Chambrion, Boscain-Chambrion-Charlot, Boscain-Charlot-2005, Boscain-Charlot-Gauthier-book, Boscain-Charlot-Gauthier, Gauthier-Jauslin, Boscain-Mason},
by Agrachev, Boscain, Chambrion, Charlot, Gauthier, Guérin, Jauslin and Mason,
\cite{Khaneja-Glaser-Brockett} by Khaneja, Glaser and Brockett,
\cite{ramakrishna-et-al-95} by Ramakrishna, Salapaka, Dahleh, Rabitz, 
\cite{Sussmann-Jurdjevic} by Sussmann and Jurdjevic, 
\cite{turinici-rabitz-01} by Turinici and Rabitz.
Let us also mention 
\cite{MM-PR-GT} by Mirrahimi, Rouchon, Turinici
and \cite{KB-JMC-MM-PR} for explicit feedback controls, inspired by Lyapunov technics.

\subsubsection{Controllability results for Schr\"odinger and wave equations}

The controllability of Schr\"odinger equations with 
distributed and boundary controls, that act linearly on the state, 
is studied since a long time.

For linear equations, the controllability is equivalent to an observability 
inequality that may be proved with different technics:
multiplier methods (see \cite{Fabre} by Fabre, \cite{Machtyngier} by Machtyngier),
microlocal analysis (see \cite{lebeau} by Lebeau, \cite{Burq} by Burq),
Carleman estimates (see \cite{Lasiecka-Triggiani, Lasiecka-Triggiani-Zhang} 
by Lasiecka, Triggiani, Zhang),
or number theory (see \cite{Tucsnak} by Ramdani, Takahashi, Tenenbaum and Tucsnak).

For nonlinear equations, we refer to
\cite{Dehman-Gerard} by Dehman, Gérard, Lebeau,
\cite{Teismann-et-alBIS} by Lange Teismann,
\cite{Laurent,Laurentdim3} by Laurent,
\cite{Rosier-Zhang} by Rosier, Zhang.

\subsubsection{Other results about bilinear quantum systems}

The study of the controllability of Schr\"odinger PDEs with 
bilinear controls started later.

The first result is negative and it is due to Turinici 
(see \cite{Defranceschi-LeBris} and Proposition \ref{Prop:BMS}).
It is a corollary of a more general result by Ball, Marsden and Slemrod 
\cite{ball-marsden-slemrod}. Because of this noncontrollability result,
such equations have been considered as non controllable for a long time.
However, important progress have been made in the last years and this
question is now better understood (see section \ref{subsub:neg}). 
Let us also mention that this negative result has been adapted to 
non linear Schrödinger equations in \cite{Teismann-et-al} by Ilner, Lange and Teismann.

Concerning exact controllability issues,
local results for 1D models have been proved in \cite{KB-JMPA,SchroLgVar} by Beauchard;
almost global results have been proved in \cite{KB-JMC}, by Coron and Beauchard.
In \cite{JMC-CRAS-Tmin}, Coron proved that a positive minimal
time was required for the local controllability of the 1D model 
(\ref{Schro_eq}) with $\mu(x)=x-1/2$.

Now, let us quote some approximate controllability results.
In \cite{KB-MM} Mirrahimi and Beauchard proved the global approximate controllability,
in infinite time, for a 1D model and in \cite{MM} Mirrahimi proved a similar result
for equations involving a continuous spectrum.
Approximate controllability, in finite time, has been proved for 
particular models by Boscain and Adami in \cite{Boscain-Adami},
by using adiabatic theory and intersection of the eigenvalues in the space of controls.
Approximate controllability, in finite time, for more general models, have been studied
by 3 teams, with different tools:
by Boscain, Chambrion, Mason, Sigalotti \cite{Chambrion-et-al},  
with geometric control methods;
by Nersesyan \cite{Nersesyan1,Nersesyan2}
with feedback controls and variational methods;
and by Ervedoza and Puel \cite{Ervedoza-Puel}
thanks to a simplified model.

Let us emphasize that the local exact controllability result of \cite{KB:SchroTIL}
and the global approximate controllability of \cite{Nersesyan1,Nersesyan2}
can be put together in order to get the
global exact controllability of 1D models (see \cite{Nersesyan2}).

Optimal control techniques have also been investigated for Schr\"odinger 
equations with a non linearity of Hartee type
in \cite{Baudouin,Baudouin-Kavian-Puel} by Baudouin, Kavian, Puel
and in \cite{CLB-Cances-Pilot} by Cances, Le Bris, Pilot.
An algorithm for the computation of such optimal controls is studied
in \cite{Baudouin-Salomon} by Baudouin and Salomon.

\subsection{Structure of this article}

This article is organized as follows.

Section \ref{sectlinearSchrod} aims at proving the controllability 
for the linear Schr\"odinger equations. 
The Subsections \ref{subsec:WP}, \ref{subsec:C1}, \ref{subsec:linearise} and 
\ref{subsect_proofthm_control_lin}
are dedicated to the different steps of the proof of Theorem \ref{Main-Thm},
where the equation is posed on a bounded interval.
The Subsection \ref{subsec:higher} is dedicated to the proof of the same result with
higher regularities, i.e. Theorem \ref{Main-ThmH5H1}.
The Subsection \ref{subsec:rad} is dedicated to the Schr\"odinger equation
for radial data on the three dimensional ball, i.e. the proof of Theorem \ref{radial-Thm}.

In Section \ref{sectNLS}, we prove Theorem \ref{thm_controlNLS} concerning the 
nonlinear Schr\"odinger equation (\ref{NLS}). 

In Section \ref{sectNLW}, we prove Theorem \ref{thm:ondes} concerning 
the nonlinear wave equation (\ref{wave}).

Finally, in Section \ref{sec:ccl}, we state some conclusions, open problems and perspectives.

\subsection{Notations}

Let us introduce some conventions and notations that are valid in all this section.
Unless otherwise specified, the functions considered are complex valued and,
for example, we write $H^1_0(0,1)$ for $H^1_0((0,1),\mathbb{C})$. When the 
functions considered are real valued, we specify it and we write, for example,
$L^2((0,T),\mathbb{R})$. We use the spaces
$$h^s(\mathbb{N}^*,\mathbb{C}):=\left\{ a=(a_k)_{k \in \mathbb{N}^*} 
\in \mathbb{C}^{\mathbb{N}^*} ;
 \sum\limits_{k=1}^{\infty} | k^s a_k |^2 < + \infty \right\}$$
equipped with the norm
$$\|a\|_{h^s}:= \Big(\sum\limits_{k=1}^{\infty} | k^s a_k |^2 \Big)^{1/2}.$$ 
The same letter $C$ denotes a positive constant, that can change from one line
to another one. If $(X,\|.\|)$ is a normed vector space and $R>0$, $B_R[X]$ denotes
the open ball $\{ x \in X ; \|x\| < R \}$ and
$\overline{B}_R[X]$ denotes the closed ball $\{ x \in X ; \|x\| \leqslant R \}$.

\section{Linear Schr\"odinger equations}
\label{sectlinearSchrod}

The goal of this section is the proof of controllability results for linear
Schr\"odinger equations, with bilinear controls. 

The Subsections \ref{subsec:WP}, \ref{subsec:C1}, \ref{subsec:linearise} and 
\ref{subsect_proofthm_control_lin}
are dedicated to the different steps of the proof of Theorem \ref{Main-Thm},
where the equation is posed on a bounded interval.
In Subsection \ref{subsec:WP}, we prove existence, uniqueness, regularity results
and bounds on the solution of the Cauchy problem (\ref{Schro_eq}), (\ref{IC}).
In Subsection \ref{subsec:C1}, we prove the $C^1$-regularity of the end-point
map associated to our control problem. 
In Subsection \ref{subsec:linearise}, we prove the controllability of the
linearized system around the ground state.
Finally, in Subsection \ref{subsect_proofthm_control_lin}, we deduce
Theorem \ref{Main-Thm} by applying the inverse mapping theorem.

The Subsection \ref{subsec:higher} is dedicated to the proof of the same result with
higher regularities, i.e. Theorem \ref{Main-ThmH5H1}.

The Subsection \ref{subsec:rad} is dedicated to the Schr\"odinger equation
for radial data on the three dimensional ball, i.e. the proof of Theorem \ref{radial-Thm}.
\\

In all this section (except in Subsection \ref{subsec:rad}), 
the operator $A$ is defined by (\ref{def:A}),
the spaces $H^s_{(0)}(0,1)$ are defined by (\ref{def:Hs}) and
$e^{-iAt}$ denotes the group of isometries of $H^s_{(0)}(0,1)$, $\forall s \geqslant 0$ 
generated by $-iA$,
\begin{equation} \label{groupe-explicit}
e^{-iAt} \varphi = \sum_{k=1}^{\infty} 
\langle \varphi , \varphi_k \rangle e^{-i \lambda_k t} \varphi_k,
\forall \varphi \in L^2(0,1).
\end{equation}
We use few classical results concerning trigonometric moment
problems that are recalled in Appendix B.

\subsection{Well posedness of the Cauchy problem}
\label{subsec:WP}

This subsection is dedicated to the statement of existence, uniqueness, 
regularity results, and bounds for the weak solutions of the Cauchy problem
\begin{equation} \label{CY}
\left\lbrace \begin{array}{l}
i \frac{\partial \psi}{\partial t} = 
- \frac{\partial^{2} \psi}{\partial x^{2}} - u(t) \mu(x) \psi - f(t,x)
, \text{  }  x \in(0,1), t \in \mathbb{R}_{+},
\\
\psi(t,0)=\psi(t,1)=0,
\\
\psi(0,x)=\psi_{0}(x).
\end{array}\right.
\end{equation}

\begin{Prop} \label{WP-CYpb}
Let $\mu \in H^{3}((0,1),\mathbb{R})$, 
$T>0$, 
$\psi_0 \in H^3_{(0)}(0,1)$, 
$f\in L^2((0,T),H^3 \cap H^1_0)$
and $u \in L^2((0,T),\mathbb{R})$.
There exists a unique weak solution of (\ref{CY}), i.e.
a function $\psi \in C^0([0,T],H^3_{(0)})$ such that the following equality 
holds in $H^{3}_{(0)}(0,1)$ for every $t \in [0,T]$,
\begin{equation} \label{SolutionFaible}
\psi(t)=e^{-iAt} \psi_{0} +i \int_{0}^{t} e^{-iA(t-\tau)}
[ u(\tau) \mu \psi(\tau) + f(\tau) ] d\tau.
\end{equation}
Moreover, for every $R>0$, there exists $C=C(T,\mu,R)>0$ such that,
if $\|u\|_{L^2(0,T)} < R$, then this weak solution satisfies
\begin{equation} \label{majo}
\| \psi \|_{C^{0}([0,T],H^{3}_{(0)})} \leqslant 
C \Big( \| \psi_{0} \|_{H^{3}_{(0)}} + \|f\|_{L^{2}((0,T),H^{3} \cap H^1_0)} \Big).
\end{equation}
If $f \equiv 0$ then 
\begin{equation} \label{NL2=1}
\|\psi(t)\|_{L^2(0,1)} = \|\psi_0\|_{L^2(0,1)}, \forall t \in [0,T].
\end{equation}
\end{Prop}

The main difficulty of the proof of this result is that
$f(s)$ is not assumed to belong to $H^3_{(0)}(0,1)$
(i.e. $f''(s,.)$ may not vanish at $x=0$ and $x=1$),
and $\mu$ is not assumed to satisfy $\mu'(0)=\mu'(1)=0$
(and thus the operator $\varphi \mapsto \mu \varphi$ does not
preserve $H^3_{(0)}(0,1)$ because for $\varphi \in H^3_{(0)}(0,1)$,
we have $(\mu \varphi)''=2 \mu' \varphi'$ at $x=0$ and $x=1$).
The argument for proving Proposition \ref{WP-CYpb} comes from the following Lemma.

\begin{Lem} \label{Lem:G}
Let $T>0$ and $f \in L^2((0,T),H^3 \cap H^1_0)$. The function
$G:t \mapsto \int_0^t e^{iAs} f(s) ds$
belongs to $C^0([0,T],H^3_{(0)})$, moreover
\begin{equation} \label{majo_Lem:G}
\|G\|_{L^\infty((0,T),H^3_{(0)})} \leqslant
c_1(T) \|f\|_{L^2((0,T),H^3 \cap H^1_0)}
\end{equation}
where the constants $c_1(T)$ are uniformly bounded for $T$ lying in bounded intervals.
\end{Lem}

\noindent \textbf{Proof of Lemma \ref{Lem:G}:} By definition, we have
$$G(t)=\sum\limits_{k=1}^{\infty} \left(
\int_0^t \langle f(s) , \varphi_k \rangle e^{i \lambda_k s} ds 
\right) \varphi_k.$$
For almost every $s \in (0,T)$, $f(s) \in H^3 \cap H^1_0$, and we have
$$\begin{array}{ll}
\langle f(s) , \varphi_k \rangle
& = \frac{1}{\lambda_k} \langle Af(s) , \varphi_k \rangle
\\ & = - \frac{\sqrt{2}}{\lambda_k} \int_0^1 f''(s,x) \sin(k\pi x) dx
\\ & = \frac{\sqrt{2}}{(k\pi)^3} \Big( (-1)^k f''(s,1)  - f''(s,0) \Big)
- \frac{\sqrt{2}}{(k\pi)^3} \int_0^1 f'''(s,x) \cos(k\pi x) dx.
\end{array}$$
Thus, we have
$$\begin{array}{lll}
\|G(t)\|_{H^3_{(0)}} 
& = &
\Big\|  \int_0^t  \langle f(s) , \varphi_k \rangle e^{i \lambda_k s} ds
\Big\|_{h^3}
\\ & \leqslant &
\frac{\sqrt{2}}{\pi^3} \left( \Big\| \int_0^t f''(s,1) e^{i\lambda_k s} ds \Big\|_{l^2}  
+ \Big\| \int_0^t f''(s,0) e^{i\lambda_k s} ds \Big\|_{l^2} 
\right)
\\ & &
+ \frac{1}{\pi^3}
\Big\| \int_0^t \langle f'''(s) ,\sqrt{2} \cos(k\pi x) \rangle e^{i\lambda_k s} 
ds \Big\|_{l^2}. 
\end{array}$$
The family $(\sqrt{2} \cos(k\pi x))_{k \in \mathbb{N}^*}$ 
is orthonormal in $L^2(0,1)$, thus 
$$\begin{array}{ll}
\Big\|  \int_0^t \langle f'''(s) , \sqrt{2}\cos(k\pi x) \rangle 
e^{i\lambda_k s} ds \Big\|_{l^2}
& =
\left( \sum\limits_{k=1}^{\infty}
\Big| \int_0^t \langle f'''(s) , \sqrt{2}\cos(k\pi x) \rangle 
e^{i\lambda_k s} ds \Big|^2 \right)^{1/2}
\\ & \leqslant
\left( \sum\limits_{k=1}^{\infty}
t  \int_0^t |\langle f'''(s) , \sqrt{2}\cos(k \pi x) \rangle|^2 ds
\right)^{1/2}
\\ & \leqslant
\sqrt{t}
\left( \int_0^t \| f'''(s)\|_{L^2}^2 ds \right)^{1/2}
\\ & \leqslant
\sqrt{t} \|f\|_{L^2((0,t),H^3)}.
\end{array}$$
Thanks to Corollary \ref{Cor3} (in Appendix B), we get
$$\begin{array}{lll}
\|G(t)\|_{H^3_{(0)}} 
& \leqslant &
\frac{\sqrt{2} C(t)}{\pi^3} \Big( \|f''(.,0)\|_{L^2(0,t)} + \|f''(.,1)\|_{L^2(0,t)} \Big) 
+ \frac{\sqrt{t}}{\pi^3} \|f\|_{L^2((0,t),H^3)}
\\ & \leqslant &
c_1(t) \|f\|_{L^2((0,t),H^3 \cap H^1_0)}
\end{array}$$
where $c_1(t)$ is uniformly bounded for $t$ lying in bounded intervals.
This bound shows that $G(t)$ belongs to $H^3_{(0)}(0,1)$ for every $t \in [0,T]$
and that the map $t\in [0,T] \mapsto G(t) \in H^3_{(0)}$ is continuous at $t=0$
(because $c_1(t)$ is uniformly bounded when $t \rightarrow 0$ and
$\|f\|_{L^2((0,t),H^3 \cap H^1_0)} \rightarrow 0$ 
when $t \rightarrow 0$, thanks to the dominated convergence theorem).
The continuity of $G$ at any $t \in (0,T)$ can be proved similarly. $\Box$
\\

\noindent \textbf{Proof of Proposition \ref{WP-CYpb}:}
Let $\mu \in H^{3}((0,1),\mathbb{R})$, $T>0$, $\psi_0 \in H^3_{(0)}(0,1)$,
$f \in L^2((0,T),H^3 \cap H^1_0)$ and $u \in L^2((0,T),\mathbb{R})$.
We consider the map
$$\begin{array}{cccc}
F: & C^0([0,T],H^3_{(0)}) & \rightarrow & C^0([0,T],H^3_{(0)}) \\
   &  \psi                & \mapsto     & \xi
\end{array}$$
where $\xi:=F(\psi)$ is defined by
\begin{equation} \label{def:xi}
\xi(t):=e^{-iAt} \psi_0 + i \int_0^t e^{-iA(t-s)} \Big(
u(s) \mu \psi(s) + f(s) \Big) ds, \forall t \in [0,T].
\end{equation}
We have assumed that $f \in L^2((0,T),H^3 \cap H^1_0)$ and $u \in L^2(0,T)$,
thus, for every $\psi \in C^0([0,T],H^3_{(0)})$, the map
$u \mu \psi + f$ belongs to
$L^2((0,T),H^3 \cap H^1_0)$ and Lemma \ref{Lem:G} ensures that
$F$ takes values in $C^0([0,T],H^3_{(0)})$. We have also used that 
in dimension $1$, $H^3$ is an algebra. 

Thanks to (\ref{majo_Lem:G}), we get, for every $t \in [0,T]$,
$$\begin{array}{ll}
\|F(\psi_1)(t)-F(\psi_2)(t)\|_{H^3_{(0)}}
& =
\Big\| \int_0^t e^{iAs} u(s) \mu (\psi_1-\psi_2)(s) ds \Big\|_{H^3_{(0)}}
\\ & \leqslant
c_1(t) \|u \mu (\psi_1-\psi_2) \|_{L^2((0,t),H^3 \cap H^1_0)}
\\ & \leqslant
c_1(t) \| u \|_{L^2(0,t)} \| \mu (\psi_1-\psi_2) \|_{L^\infty((0,t),H^3 \cap H^1_0)}
\\ & \leqslant
c_1(t) \|u\|_{L^2(0,t)} C(\mu) \| \psi_1 - \psi_2 \|_{L^\infty((0,t),H^3_{(0)})}
\end{array}$$
thus
\begin{equation} \label{def:c_2}
\| F(\psi_1)-F(\psi_2) \|_{L^\infty((0,T),H^3_{(0)})}
\leqslant c_2(T,\mu) 
\|u\|_{L^2(0,T)} \| \psi_1 - \psi_2 \|_{L^\infty((0,T),H^3_{(0)})}.
\end{equation}
If $\|u\|_{L^2(0,T)}$ is small enough, then
$F$ is a contraction. Thanks to the Banach fixed point theorem,
there exists $\psi \in C^0([0,T],H^3_{(0)})$ such that $F(\psi)=\psi$.
The previous arguments show that, for this fixed point, we have
$$\| \psi \|_{L^\infty((0,T),H^3_{(0)})}
 \leqslant \|\psi_0\|_{H^3_{(0)}} + 
c_2(T,\mu) \|u\|_{L^2(0,T)} \| \psi \|_{L^\infty((0,T),H^3_{(0)})} +
c_1(T) \|f\|_{L^2((0,T),H^3 \cap H^1_0)}.$$
Thus, if $c_2(T,\mu) \|u\|_{L^2(0,T)} \leqslant 1/2$,
then, we get (\ref{majo}). 

We have proved Proposition \ref{WP-CYpb} when $\|u\|_{L^2(0,T)}$ is small enough.
If it is not the case, one may consider $0=T_0<T_1<...<T_N=T$ such that
$\|u\|_{L^2(T_j,T_{j+1})}$ is small and apply the previous result on
$[T_0,T_1],...,[T_{N-1},T_N]$ in order to get the conclusion. 
Since our constant $c_1(t)$ is uniform on bounded sets, 
we easily get that $N$ only depends on $R$, so that the constant in 
the Proposition does only depend on $T$, $\mu$ and $R$ as claimed.

Now, let us prove that (\ref{NL2=1}) holds when $f=0$.
Classical arguments allow to prove that, when $u \in C^0([0,T],\mathbb{R})$,
then $\psi \in C^1([0,T],L^2)$ and the first equality of (\ref{Schro_eq})
holds in $L^2$ for every $t \in [0,T]$. Thus, when $u \in C^0([0,T],\mathbb{R})$,
we can take the $L^2$-scalar product of this equation with $\psi$;
and the imaginary part of the resulting equality gives
$$\frac{d}{dt} \| \psi(t) \|_{L^2}^2 =0.$$
Thus, we have (\ref{NL2=1}) when $u \in C^0([0,T],\mathbb{R})$.
A density argument allows to prove (\ref{NL2=1}) 
when $u$ only belongs to $L^2((0,T),\mathbb{R})$. $\Box$

\subsection{$C^1$-regularity of the end-point map}
\label{subsec:C1}

For $T>0$ we introduce the tangent space of $\mathcal{S}$ at $\psi_1(T)$
$$V_T:=\{ \xi \in L^2(0,1) ; \Re \langle \xi , \psi_1(T) \rangle =0 \}$$
and the orthogonal projection
$$P_T : L^2(0,1) \rightarrow V_T.$$
Proposition \ref{WP-CYpb} allows to consider the map
\begin{equation} \label{def:ThetaT}
\begin{array}{cccc}
\Theta_T : & L^2((0,T),\mathbb{R}) & \rightarrow & V_T \cap H^3_{(0)}(0,1) \\
           &               u       & \mapsto     & P_T [\psi(T)]
\end{array}
\end{equation}
where $\psi$ is the solution of (\ref{Schro_eq}), (\ref{IC}).
The goal of this section is the proof of the following result.

\begin{Prop} \label{Prop:C1}
Let $T>0$ and $\mu \in H^{3}((0,1),\mathbb{R})$.
The map $\Theta_T$ defined by (\ref{def:ThetaT}) is $C^1$.
Moreover, for every $u, v \in L^2((0,T),\mathbb{R})$, we have
\begin{equation} \label{dThetaT}
d\Theta_T(u).v=P_T[\Psi(T)]
\end{equation}
where $\Psi$ is the weak solution of the linearized system
\begin{equation} \label{proofC1:L}
\left\lbrace\begin{array}{l}
i \frac{\partial \Psi}{\partial t} = - \Psi'' - u(t) \mu(x) \Psi
- v(t) \mu(x) \psi , x \in (0,1), t \in (0,T),
\\
\Psi(t,0)=\Psi(t,1)=0,
\\
\Psi(0,x)=0,
\end{array}\right.
\end{equation} 
and $\psi$ is the solution of (\ref{Schro_eq}),(\ref{IC}).
\end{Prop}

\noindent \textbf{Proof of Proposition \ref{Prop:C1}:}
Let $T>0$, $\mu \in H^{3}((0,1),\mathbb{R})$ and $u \in L^2((0,T),\mathbb{R})$.
First, let us emphasize that the linear map $v \mapsto \Psi(T)$ is
continuous from $L^2((0,T),\mathbb{R})$ to $H^3_{(0)}(0,1)$ 
thanks to Proposition \ref{WP-CYpb}.

\emph{First step: We prove that $\Theta_T$ is differentiable and 
that (\ref{dThetaT}) holds.} Let $\psi$ be the weak solution
of (\ref{Schro_eq}),(\ref{IC}), $\Psi$ solution of (\ref{proofC1:L}) and $\widetilde{\psi}$ solution of
\begin{equation} \label{eq:psi-tilde}
\left\lbrace\begin{array}{l}
i \frac{\partial \widetilde{\psi}}{\partial t} = 
- \widetilde{\psi}'' - (u+v)(t) \mu(x) \widetilde{\psi}, x \in (0,1), t \in (0,T),
\\
\widetilde{\psi}(t,0)=\widetilde{\psi}(t,1)=0,
\\
\widetilde{\psi}(0,x)=\varphi_1.
\end{array}\right.
\end{equation}
Then $\Delta:=\widetilde{\psi}-\psi-\Psi$ is the weak solution of
\begin{equation} \label{eq:Delta}
\left\lbrace\begin{array}{l}
i \frac{\partial \Delta}{\partial t} = - \Delta'' - (u+v)(t) \mu(x) \Delta - v(t) \mu \Psi
, x \in (0,1), t \in (0,T),
\\
\Delta(t,0)=\Delta(t,1)=0,
\\
\Delta(0,x)=0.
\end{array}\right.
\end{equation}
Let us prove that
\begin{equation} \label{proofC1:step1}
\| \Delta \|_{C^0([0,T],H^3_{(0)})} = o ( \|v\|_{L^2} )
\text{ when } \|v\|_{L^2} \rightarrow 0,
\end{equation}
which gives the conclusion. Let $R>0$ be such that $\|u\|_{L^2(0,T)} < R$
and $\|u+v\|_{L^2(0,T)}<R$. Thanks to Proposition \ref{WP-CYpb}, 
there exists $C_j=C_j(T,\mu,R)>0$ for $j=0,1$ such that
$$\| \Delta \|_{C^0([0,T],H^3_{(0)})} \leqslant
C_0 \| v \mu \Psi \|_{L^2((0,T),H^3 \cap H^1_0)} \leqslant
C_1 \|v\|_{L^2} \|\Psi \|_{C^0([0,T],H^3_{(0)})},$$
$$\begin{array}{ll}
\| \Psi \|_{C^0([0,T],H^3_{(0)})} 
& \leqslant C_0 \| v \mu \psi \|_{L^2((0,T),H^3 \cap H^1_0)}
\\ 
& \leqslant C_1 \|v\|_{L^2} \| \psi \|_{C^0([0,T],H^3_{(0)})} 
\\ 
& \leqslant C_0 C_1 \|v\|_{L^2} \|\varphi_1\|_{H^3_{(0)}},
\end{array}$$
which proves (\ref{proofC1:step1}).
\\

\emph{Second step: We prove that $d\Theta_T$ is continuous.}
Actually, we prove that this map is locally Lipschitz.
Let $u, \widetilde{u} \in L^2((0,T),\mathbb{R})$ and 
$v \in L^2((0,T),\mathbb{R})$.
Let $\psi$ be the solution of (\ref{Schro_eq}),(\ref{IC}), $ \Psi$ solution of (\ref{proofC1:L}) and $\widetilde{\psi}$, $\widetilde{\Psi}$  solution of 
$$\begin{array}{ll}
\left\lbrace \begin{array}{l}
i \frac{\partial \widetilde{\psi}}{\partial t} = - \widetilde{\psi}'' - 
\widetilde{u}(t) \mu(x) \widetilde{\psi}, 
\\
\widetilde{\psi}(t,0)=\widetilde{\psi}(t,1)=0,
\\
\widetilde{\psi}(0,x)=\varphi_1,
\end{array}\right.
&
\left\lbrace \begin{array}{l}
i \frac{\partial \widetilde{\Psi}}{\partial t} = 
- \widetilde{\Psi}'' - \widetilde{u}(t) \mu(x) \widetilde{\Psi}
- v(t) \mu(x) \widetilde{\psi} ,
\\
\widetilde{\Psi}(t,0)=\widetilde{\Psi}(t,1)=0,
\\
\widetilde{\Psi}(0,x)=0,
\end{array}\right.
\end{array}$$
We have 
$$[d\Theta_T(u)-d\Theta_T(\widetilde{u})].v=P_T[ \Psi(T)-\widetilde{\Psi}(T)]
= P_T[ \Xi(T) ]$$
where $\Xi$ is the weak solution of
$$\left\lbrace \begin{array}{l}
i \frac{\partial \Xi}{\partial t}=
- \frac{\partial^2 \Xi}{\partial x^2}
- u(t) \mu \Xi
- (u-\widetilde{u}) \mu \widetilde{\Psi}
- v \mu (\psi-\widetilde{\psi}),\\
\Xi(t,0)=\Xi(t,1)=0,\\
\Xi(0)=0.
\end{array}\right.$$
Let $R>0$ be such that
$\|u\|_{L^2(0,T)}<R$, $\| \widetilde{u} \|_{L^2(0,T)}<R$.
Let us prove that
$$\| \Xi  \|_{C^0([0,T],H^3_{(0)})}
\leqslant \mathcal{C} \|v\|_{L^2} \|u-\widetilde{u}\|_{L^2}$$
where $\mathcal{C}=\mathcal{C}(T,\mu,R)>0$, which gives the conclusion.
Thanks to Proposition \ref{WP-CYpb}, we have
$$\begin{array}{ll}
\| \Xi  \|_{C^0([0,T],H^3_{(0)})}
& \leqslant C_2 
\| (u-\widetilde{u}) \mu \widetilde{\Psi} 
+  v \mu (\psi-\widetilde{\psi}) \|_{L^2((0,T),H^3 \cap H^1_0)}
\\ & \leqslant
C_3 \Big( 
\|u-\widetilde{u}\|_{L^2} \| \widetilde{\Psi} \|_{C^0([0,T],H^3_{(0)})} 
+ \|v\|_{L^2} \| \psi-\widetilde{\psi} \|_{C^0([0,T],H^3_{(0)})} 
\Big)
\\ & \leqslant
C_4 \Big( 
\|u-\widetilde{u}\|_{L^2} \| v \mu \widetilde{\psi} \|_{L^2((0,T),H^3 \cap H^1_0)}
+ \|v\|_{L^2} \| (\widetilde{u}-u) \mu \widetilde{\psi} \|_{L^2((0,T),H^3 \cap H^1_0)}
\Big)
\\ & \leqslant
C_5 \Big(
\|u-\widetilde{u}\|_{L^2} \| v \|_{L^2} \| \widetilde{\psi} \|_{C^0([0,T],H^3_{(0)})}+
\|v\|_{L^2} \| \widetilde{u}-u \|_{L^2} \| \widetilde{\psi} \|_{C^0([0,T],H^3_{(0)})}
\Big)
\\ & \leqslant
C_6 \|u-\widetilde{u}\|_{L^2} \| v \|_{L^2},
\end{array}$$
where $C_j=C_j(T,\mu,R)>0$ for $j=2,...,6$. $\Box$

\subsection{Controllability of the linearized system}
\label{subsec:linearise}

The goal of this section is the proof of the following result.

\begin{Prop} \label{Cont-Lin-H3L2}
Let $T>0$ and $\mu \in H^{3}((0,1),\mathbb{R})$ be such that (\ref{hyp_mu}) holds.
The linear map $d\Theta_T(0):L^2((0,T),\mathbb{R}) \rightarrow V_T \cap H^3_{(0)}(0,1)$
has a continuous right inverse
$d\Theta_T(0)^{-1}:V_T \cap H^3_{(0)}(0,1) \rightarrow L^2((0,T),\mathbb{R})$.
\end{Prop}

The proof of Proposition \ref{Cont-Lin-H3L2} relies on an 
Ingham inequality, due to Haraux (see \cite{Haraux} and Appendix B).
\\

\noindent \textbf{Proof of Proposition \ref{Cont-Lin-H3L2}:}
We have $d\Theta_T(0).v=\Psi(T)$ where
\begin{equation} \label{syst-lin}
\left\lbrace\begin{array}{l}
i \frac{\partial \Psi}{\partial t} = - \Psi'' - v(t) \mu \psi_{1} , 
\\
\Psi(t,0)=\Psi(t,1)=0,
\\
\Psi(0)=0,
\end{array}\right.
\end{equation}
thus
$$\Psi(T)=\sum_{k=1}^{\infty}
i \langle \mu \varphi_1 , \varphi_k \rangle 
\left(\int_0^T v(t) e^{i(\lambda_k-\lambda_1)t} dt\right) e^{-i\lambda_k T} \varphi_k.$$
Let $\Psi_{f} \in V_T \cap H^{3}_{(0)}(0,1)$.
If $\Psi$ is the solution of (\ref{syst-lin}) for some
$v \in L^2((0,T),\mathbb{R})$, then, the equality $\Psi(T)=\Psi_{f}$ 
is equivalent to the trigonometric moment problem
\begin{equation} \label{CL:trigo-moment-pb}
\int_{0}^{T} v(t) e^{i(\lambda_{k}-\lambda_{1})t} dt
= d_{k-1}(\Psi_f):=
\frac{\langle \Psi_{f} , \varphi_{k} \rangle e^{i \lambda_{k}T}}{
i \langle \mu \varphi_{1} , \varphi_{k} \rangle},
 \forall k \in \mathbb{N}^{*}.
\end{equation}
Now, we apply Corollary \ref{Cor:haraux1} with
$\omega_k := \lambda_{k+1}-\lambda_1, \forall k \in \mathbb{N}$,
and we get the conclusion with
$$d\Theta_T(0)^{-1}(\Psi_f):=L[d(\Psi_f)],$$
where $d(\Psi_f):=( d_{k}(\Psi_f) )_{k \in \mathbb{N}}$.
Indeed, for $\Psi_f \in V_T \cap H^3_{(0)}(0,1)$, the sequence
$d(\Psi_f)$ belongs to $l^2_r(\mathbb{N},\mathbb{C})$ thanks to
the assumption (\ref{hyp_mu}). $\Box$

\subsection{Proof of Theorem \ref{Main-Thm}}
\label{subsect_proofthm_control_lin}

Let $T>0$ and $\mu \in H^{3}((0,1),\mathbb{R})$ be such that 
(\ref{hyp_mu}) holds. 
Let $R_1 >0$ and $\delta>0$ be such that,
$$ \forall u \in B_{R_1}[L^2((0,T),\mathbb{R})], \text{ the solution of } 
(\ref{Schro_eq}), (\ref{IC}) \text{ satifies } 
\Re \langle \psi(T) , \psi_1(T) \rangle >0,$$
(see Proposition \ref{WP-CYpb}) and
$$ \forall \psi_f \in \mathcal{S} \cap H^3_{(0)}(0,1)
\text{ with } \| \psi_f - \psi_1(T)\|_{H^3_{(0)}}<\delta,
\text{ we have } \Re \langle \psi_f , \psi_1(T) \rangle >0.$$ 
The spaces $\overline{B}_{R_1}[L^2((0,T),\mathbb{R})]$ and 
$V_T \cap H^3_{(0)}(0,1)$ are Banach spaces.
The map $\Theta_T: \overline{B}_{R_1}[L^2((0,T),\mathbb{R})] 
\rightarrow V_T \cap H^3_{(0)}(0,1)$ is $C^1$
(see Proposition \ref{Prop:C1}),
its differential at $0$ has a continuous right inverse
$d\Theta_T(0)^{-1}:V_T \cap H^3_{(0)}(0,1) \rightarrow L^2((0,T),\mathbb{R})$
(see Proposition \ref{Cont-Lin-H3L2}). 
Thanks to the inverse mapping theorem,
there exists $\delta_1 \in (0,\delta)$ and a $C^1$ map
$$\Theta_T^{-1}:B_{\delta_1}[ V_T \cap H^3_{(0)}(0,1)] \rightarrow 
\overline{B}_{R_1}[L^2((0,T),\mathbb{R})]$$
such that $\Theta_T( \Theta_T^{-1} (\widetilde{\psi_f} ))=\widetilde{\psi_f}$
for every $\widetilde{\psi}_f \in B_{\delta_1}[ V_T \cap H^3_{(0)}(0,1)]$.

For $\psi_f \in \mathcal{S} \cap H^3_{(0)}(0,1)$ 
with $\| \psi_f - \psi_1(T)\|_{H^3_{(0)}}<\delta_1$,
we have $\| P_T \psi_f \|_{H^3_{(0)}} < \delta_1$, 
thus we can define
$$\Gamma(\psi_f)=:\Theta_T^{-1}[P_T \psi_f ].$$
Thanks to the choice of $R_1$ and $\delta$ we know that the solution of
(\ref{Schro_eq}), (\ref{IC}) with $u=\Gamma(\psi_f)$ satisfies
$$\begin{array}{ll}
\psi(T)
& = P_T ( \psi(T) ) + \sqrt{1-\| P_T \psi(T) \|_{L^2}^2} \psi_1(T) \\
& = P_T(\psi_f) + \sqrt{1-\| P_T \psi_f \|_{L^2}^2} \psi_1(T) 
= \psi_f.
\end{array}$$

\subsection{Generalization to higher regularities}
\label{subsec:higher}

The goal of this section is the proof of Theorem \ref{Main-ThmH5H1}.
The first step of the proof consists in adapting Proposition \ref{WP-CYpb}.

\begin{Prop} \label{WP-CYpbH5H1}
Let $\mu \in H^{5}((0,1),\mathbb{R})$, 
$T>0$, 
$\psi_0 \in H^5_{(0)}(0,1)$, 
$f\in H^1_0((0,T),H^3 \cap H^1_0)$
and $u \in H^1_0((0,T),\mathbb{R})$.
There exists a unique function
$\psi \in C^1([0,T],H^3_{(0)})$ 
such that the equality (\ref{SolutionFaible})
holds in $C^1([0,T],H^{3}_{(0)})$.
Moreover, for every $R>0$ there exists $C=C(T,\mu,R)>0$ such that,
if $\|u\|_{H^1_0(0,T)}<R$, then, this weak solution satisfies
\begin{equation} \label{borne-H2}
\| \psi \|_{C^{1}([0,T],H^{3}_{(0)})} \leqslant 
C \Big( \| \psi_{0} \|_{H^{5}_{(0)}} + \|f\|_{H^{1}((0,T),H^{3} \cap H^1_0)} \Big).
\end{equation}
\end{Prop}

The proof of Proposition \ref{WP-CYpbH5H1} is the same as the one of
Proposition \ref{WP-CYpb}, except that we use the following Lemma, 
instead of Lemma \ref{Lem:G}.

\begin{Lem} \label{Lem:G-H5}
Let $T>0$, $u_0\in H^5 \cap H_{(0)}^3$ and $f \in H^1((0,T),H^3 \cap H^1_0)$ be such that $-iA u_0+ f(0)\in H_{(0)}^3$.
The function $G:t \mapsto e^{-iAt}u_0+\int_0^t e^{-iA(t-s)} f(s) ds$
belongs to $C^1([0,T],H^3_{(0)})$, moreover
$$\|G\|_{C^1([0,T],H^3_{(0)})} \leqslant
c_1(T) \left(\|u_0\|_{H^3_{(0)}}+\|f\|_{H^1((0,T),H^3 \cap H^1_0)}+\|-iA u_0+ f(0)\|_{H^3_{(0)}} \right)$$
where the constants $c_1(T)$ are uniformly bounded for $T$ lying in bounded intervals. 
We also have
$$\|-iA G(T)+ f(T)\|_{H^3_{(0)}} \leqslant c_1(T) \left(\|u_0\|_{H^3_{(0)}}+\|f\|_{H^1((0,T),H^3 \cap H^1_0)}+\|-iA u_0+ f(0)\|_{H^3_{(0)}} \right).$$
\end{Lem}

\noindent \textbf{Proof of Lemma \ref{Lem:G-H5}:}
We already know that $G \in C^0([0,T],H^3_{(0)})$.
First let us write
$$G(t)=e^{-iAt}u_0+\int_0^t e^{-iA\tau} f(t-\tau) d\tau.$$
Since $u_0 \in H^4_{(0)}$ and $f \in H^1((0,T),H^2_{(0)})$,
we know that $G \in C^1([0,T],H^2_{(0)})$ and
the following equality holds in $H^2_{(0)}$ for every $t \in [0,T]$,
\begin{eqnarray*}
\frac{\partial G}{\partial t}(t)&=&-iAe^{-iAt}u_0+e^{-iAt}f(0)+
\int_0^t e^{-iA\tau} \frac{\partial f}{\partial t}(t-\tau) d\tau\\
&=&e^{-iAt}\left[-iAu_0+f(0)\right]+ \int_0^t e^{-iA(t-s)} 
\frac{\partial f}{\partial t}(s) ds
\end{eqnarray*}
(the proof of this result involves classical technics). 
Thanks to this expression and Lemma \ref{Lem:G}, we get 
$$\frac{\partial G}{\partial t} \in C^0([0,T],H^3_{(0)}).$$

Let us prove that $G \in C^1([0,T],H^3_{(0)})$, i.e.
for every $t \in [0,T]$,
$$\Big\| \frac{G(t+h)-G(t)}{h} - \frac{\partial G}{\partial t}(t)
\Big\|_{H^3_{(0)}} \rightarrow 0
\text{ when } h \rightarrow 0.$$
We have
\begin{equation} \label{diff-H3}
\begin{array}{ll}
\frac{G(t+h)-G(t)}{h} - \frac{\partial G}{\partial t}(t)
= & 
e^{-iAt} \left[ \frac{e^{-iAh}-Id}{h} u_0 + i Au_0 - f(0) \right]
+ \frac{1}{h} \int_{t}^{t+h} e^{-iA\tau} f(t+h-\tau) d\tau
\\ &
+ \int_0^t e^{-iA\tau} \left[
\frac{f(t+h-\tau)-f(t-\tau)}{h} - \frac{\partial f}{\partial t}(t-\tau)
\right] d\tau.
\end{array}
\end{equation}
By applying Lemma \ref{Lem:G}, we see that
the $H^3_{(0)}(0,1)$-norm of the second term of the right hand side of (\ref{diff-H3})
tends to zero when $h \rightarrow 0$ because $f \in H^1((0,T),H^3 \cap H^1_0)$.
Thanks to several changes of variables, 
the first term of the right hand side of (\ref{diff-H3})
may be decomposed in the following way 
\begin{equation} \label{diff-H3-BIS}
\begin{array}{l}
e^{-iAt} \left[ \frac{e^{-iAh}-Id}{h} \Big( u_0 +i A^{-1} f(0) \Big)
+ i A \Big( u_0 +i A^{-1} f(0) \Big) \right] 
\\
+ e^{-iAt} \frac{1}{h} \int_0^h e^{-iAs} \Big( f(h-s) - f(0) \Big) ds.
\end{array}
\end{equation}
The $H^3_{(0)}(0,1)$-norm of the first term of (\ref{diff-H3-BIS})
tends to zero when $h \rightarrow 0$ because
$u_0 +i A^{-1} f(0) \in H^5_{(0)}(0,1)$.
The $H^3_{(0)}(0,1)$-norm of the second term of (\ref{diff-H3-BIS})
also tends to zero when $h \rightarrow 0$ because,
thanks to Lemma \ref{Lem:G} and Cauchy-Schwarz inequality, it is bounded by
$$\begin{array}{ll}
& \Big\| 
\int_0^h e^{iAs} \Big( 
\frac{f(s) - f(0)}{h}\Big) ds
\Big\|_{H^3_{(0)}}\leqslant 
c_0(h) \Big\|  
\frac{f(.) - f(0)}{h} 
\Big\|_{L^2((0,h),H^3 \cap H^1_0)}
\\ 
\leqslant &
\frac{c_0(h)}{h}\Big\| \int_0^{.} \frac{\partial f}{\partial t}(\tau)d\tau 
\Big\|_{L^{2}((0,h),H^3 \cap H^1_0)} \leqslant \frac{c_0(h)\sqrt{h}}{h}\Big\| \int_0^{.} \frac{\partial f}{\partial t}(\tau)d\tau 
\Big\|_{L^{\infty}((0,h),H^3 \cap H^1_0)} \\
\leqslant & c_0(h) \Big\| \frac{\partial f}{\partial t}\Big\|_{L^{2}((0,h),H^3 \cap H^1_0)}. 

\end{array}$$

The estimate (\ref{majo_Lem:G}) of Lemma \ref{Lem:G} gives the first inequality
of Lemma \ref{Lem:G-H5}. 
Moreover, by integration by part in time, we get
\begin{eqnarray*}
-iAG(t)&=&-iAe^{-iAt}u_0-\int_0^t iAe^{-iA\tau} f(t-\tau) d\tau\\
&=&iAe^{-iAt}u_0+e^{-iAt}f(0)-f(t)+\int_0^t e^{-iA\tau} \frac{\partial f}{\partial t}(t-\tau) d\tau.
\end{eqnarray*}
and we get the second estimate thanks to the identity
\begin{eqnarray*}
-iAG(t)+f(t)=e^{-iAt}\left[iAu_0+f(0)\right]+\int_0^t e^{-iA(t-\tau)} \frac{\partial f}{\partial t}(\tau) d\tau. \Box
\end{eqnarray*}

The following statement is the appropriate adaptation of Propositon \ref{Prop:C1}.

\begin{Prop} \label{Prop:C1-H5}
Let $T>0$ and $\mu \in H^{5}((0,1),\mathbb{R})$.
The map $\Theta_T$ defined by (\ref{def:ThetaT}) is $C^1$
from $H^1_0((0,T),\mathbb{R})$ to $V_T \cap H^5_{(0)}(0,1)$.
\end{Prop}

\noindent \textbf{Proof of Proposition \ref{Prop:C1-H5}:}

\emph{First step: we prove that $\Theta_T$ maps $H^1_0((0,T),\mathbb{R})$ into 
$V_T \cap H^5_{(0)}(0,1)$.} Let $u \in H^1_0((0,T),\mathbb{R})$ and $\psi$ be the
weak solution of (\ref{Schro_eq}), (\ref{IC}). Then
$\psi \in C^1([0,T],H^2_{(0)}) \cap C^0([0,T],H^4_{(0)})$ and
the first equality of (\ref{Schro_eq}) holds in $H^2_{(0)}$ for every $t \in [0,T]$
(the proof of this result involves classical technics).
In particular, we have
$$\begin{array}{ll}
\| \psi(T) \|_{H^5_{(0)}}
& = \| \psi''(T)\|_{H^3_{(0)}} \\
& = \Big\| \frac{\partial \psi}{\partial t}(T) \Big\|_{H^3_{(0)}}
\text{ because } u(T)=0
\end{array}$$
which is finite, thanks to Proposition \ref{WP-CYpbH5H1}.
Thus, $\Theta_T$ maps $H^1_0((0,T),\mathbb{R})$ into $V_T \cap H^5_{(0)}(0,1)$.
\\

\emph{Second step: We prove that 
$\Theta_T:H^1_0((0,T),\mathbb{R}) \rightarrow V_T \cap H^5_{(0)}$ is differentiable.}
Let $u,v \in H^1_0((0,T),\mathbb{R})$, $\psi$, $\Psi$, 
$\widetilde{\psi}$ be the weak solutions of
(\ref{Schro_eq}),(\ref{IC}), (\ref{proofC1:L}), (\ref{eq:psi-tilde}).
Then, $\Delta:=\widetilde{\psi}-\psi-\Psi$ is the weak solution of (\ref{eq:Delta}).
Let us prove that
$$\| \Delta(T) \|_{H^5_{(0)}} = o ( \|v\|_{H^1_0} )
\text{ when } \|v\|_{H^1_0} \rightarrow 0,$$
which gives the conclusion. Let $R>0$ be such that
$\|u\|_{H^1_0}<R$ and $\|u+v\|_{H^1_0}<R$. Thanks to Proposition \ref{WP-CYpbH5H1},
there exists $C=C(T,\mu,R)>0$, $C_1=C_1(\mu)>0$ such that
$$\begin{array}{ll}
\| \Delta(T)\|_{H^5_{(0)}}
& = \| \Delta''(T) \|_{H^3_{(0)}} \\
& = \Big\|  \frac{\partial \Delta}{\partial t}(T) \Big\|_{H^3_{(0)}}
\text{ because } u(T)=v(T)=0 \\
& \leqslant C \| v \mu \Psi \|_{H^1_0((0,T),H^3 \cap H^1_0)} \\
& \leqslant C C_1 \|v\|_{H^1_0} \|\Psi\|_{C^1([0,T],H^3_{(0)})} \\
& \leqslant C^2 C_1 \|v\|_{H^1_0} \| v \mu \psi \|_{H^1_0((0,T),H^3 \cap H^1_0)} \\
& \leqslant C^2 C_1^2 \|v\|_{H^1_0}^2 \| \psi \|_{C^1([0,T],H^3_{(0)})}.
\end{array}$$
The proof of the continuity of the map 
$d\Theta_T: H^1_0((0,T),\mathbb{R}) \rightarrow \mathcal{L}( H^1_0 , V_T \cap H^5_{(0)})$
involves similar arguments. $\Box$

\begin{rk}
With the same kind of arguments, we could get that 
$A\psi +u(t)\mu \psi \in C^0([0,T],H^3_{(0)})$. Therefore,  $\psi(t)$ does not, in general, belong to $H^5_{(0)}(0,1)$ 
for $t\in (0,T)$. 
\end{rk}

The following statement is the appropriate generalization of Proposition \ref{Cont-Lin-H3L2}.

\begin{Prop} \label{Control-Linearise}
Let  $T>0$, $\mu \in H^{5}((0,1),\mathbb{R})$ be such that
(\ref{hyp_mu}) holds and $\Theta_T$ be defined by (\ref{def:ThetaT}).
The linear map $d\Theta_T(0):H^1_0((0,T),\mathbb{R}) \rightarrow V_T \cap H^5_{(0)}(0,1)$ 
has a continuous right inverse
$d\Theta_T(0)^{-1} : V_T \cap H^5_{(0)}(0,1) \rightarrow  H^1_0((0,T),\mathbb{R})$.
\end{Prop}

\noindent \textbf{Proof of Proposition \ref{Control-Linearise}:}
Let $\Psi_f \in V_T \cap H^5_{(0)}(0,1)$.
If $\Psi$ is the solution of (\ref{syst-lin}) for some
$v \in H^1_0((0,T),\mathbb{R}$, then, the equality $\Psi(T)=\Psi_{f}$ 
is equivalent to the trigonometric moment problem (\ref{CL:trigo-moment-pb}), 
or equivalently
\begin{equation} 
\begin{array}{l}
\int_0^T \dot{v}(t) dt = 0,
\\
\int_0^T (T-t) \dot{v}(t)dt= 
\frac{1}{i\langle \mu \varphi_{1} , \varphi_{1} \rangle}
\langle \Psi_{f} , \varphi_{1} \rangle e^{i \lambda_{1}T},
\\
\int_{0}^{T} \dot{v}(t) e^{i(\lambda_{k}-\lambda_{1})t} dt
= \frac{\lambda_{1}-\lambda_{k}}{ \langle \mu \varphi_{1} , \varphi_{k} \rangle}
\langle \Psi_{f} , \varphi_{k} \rangle e^{i \lambda_{k}T}, \forall k \geqslant 2.
\end{array}
\end{equation}
The conclusion comes from Corollary \ref{Cor:haraux} (in Appendix B). $\Box$
\\

Now, Theorem \ref{Main-ThmH5H1} may be proved exactly as Theorem \ref{Main-Thm}.

\subsection{Case of the three dimensional ball with radial data}
\label{subsec:rad}

The goal of this section is the proof of Theorem \ref{radial-Thm}.
This proof is very similar to the case of the interval 
and we only give the necessary modifications. The equivalent of 
Lemma \ref{Lem:G} is proved with a similar computation for 
$f\in L^2((0,T),H^3_{rad}\cap H^1_{(0)})$. More precisely, 
for almost every $s\in (0,T)$, we have
\begin{eqnarray*}
\left\langle f(s),\varphi_k\right\rangle &=&\int_{B^3} f(s)\varphi_k =\frac{1}{\lambda_k^2}\int_{B^3} f(s) \Delta^2 \varphi_k=\frac{1}{\lambda_k^2}\int_{B^3} \Delta f(s)  \Delta \varphi_k\\
&=&-\frac{1}{\lambda_k^2}\int_{B^3} \nabla \Delta f(s)\cdot \nabla \varphi_k +\frac{1}{\lambda_k^2}\int_{S^2} \Delta f(s) \frac{\partial \varphi_k }{\partial n} d\sigma.
\end{eqnarray*}
To bound the first term, we use $\nabla \Delta f\in L^2((0,T),L^2(B^3)^3)$ and the fact that the functions 
$( \nabla \varphi_k / \sqrt{\lambda_k} )_{k\in \mathbb{N}^*}$
form an orthonormal family of $L^2(B^3)^3$ because
$$
\int_{B^3} \nabla \varphi_i\cdot  \nabla \varphi_j = - \int_{B^3} \varphi_i \Delta\varphi_j= \lambda_j \delta_{i,j}.
$$
For the second term, since $f$ and $\varphi_k$ are radial, we have 
$$
\frac{1}{\lambda_k^2}\int_{S^2} \Delta f(s) \frac{\partial \varphi_k }{\partial n} d\sigma =\frac{2^{3/2}\sqrt{\pi}(-1)^k}{\lambda_k^{3/2}} \Delta f(s,r=1).
$$
We conclude as in Lemma \ref{Lem:G} for this term since the eigenvalues are the same and Corollary \ref{Cor3} still applies. The genericity of assumption (\ref{hyp_mu_rad}) is detailed in the Appendix A, Proposition \ref{prop_hyp_mu-generic}.
\begin{rk}
It is very likely that the same analysis would work in any dimension $n\leqslant 5$, provided that $H^3$ remains an algebra. However, this would require the analysis of the zeros of the Bessel functions and we have chosen to present the simplest result.
\end{rk}

\section{Nonlinear Schr\"odinger equations}
\label{sectNLS}

In this section, we study the nonlinear Schr\"odinger equation
with Neumann boundary conditions (\ref{NLS}). 
The goal is the proof of Theorem \ref{thm_controlNLS}
\\

First, let us introduce the following notations, 
that will be valid in all the section \ref{sectNLS}. The operator $A$ is defined by
\begin{equation} \label{def:A-NLS}
\begin{array}{ll}
D(A)=H^2_{(0)}(0,1):=\{ \varphi \in H^2(0,1) ; \varphi'(0)=\varphi'(1)=0 \}, &
A \varphi = - \varphi''.
\end{array}
\end{equation}
Its eigenvectors $(\varphi_{k})_{k \in \mathbb{N}}$ 
and eigenvalues $(\lambda_{k})_{k \in \mathbb{N}}$ are
\begin{equation} \label{vep-NLS}
\begin{array}{ll}
\varphi_{0} := 1 , & \lambda_{0}:=0 \\
\varphi_k(x) :=\sqrt{2} \cos(k\pi x) & \lambda_{k}:=(k\pi)^2, \forall k \in \mathbb{N}^*.
\end{array}
\end{equation}
We introduce the spaces
\begin{equation} \label{def:Hs-NLS}
H^s_{(0)}(0,1):=D( A^{s/2} ), \forall s>0
\end{equation}
and the notation 
\begin{equation} \label{def:k*}
k_*:=\max\{ k , 1 \}, \forall k \in \mathbb{N}.
\end{equation}

\subsection{Well posedness of the Cauchy problem}

The goal of this subsection is the proof of the following result.

\begin{Prop} \label{WP-NLS}
Let $\mu \in H^2((0,1),\mathbb{R})$ and $T>0$.
There exists $\delta>0$ such that,
for every $u \in B_\delta [L^2(0,T)]$,
there exists a unique weak solution $\psi \in C^0([0,T],H^2_{(0)})$ 
of (\ref{NLS}), (\ref{IC-NLS}).
Moreover, we have
\begin{equation*}
\|\psi(t)\|_{L^2(0,1)} = \|\psi_0\|_{L^2(0,1)}, \forall t \in [0,T].
\end{equation*}
\end{Prop}

We search $\psi$ in the form $\psi(t,x)=e^{-it}(1+\zeta(t,x))$,
where $\zeta$ is a weak solution of 
\begin{equation} \label{eq:zeta}
\left\lbrace \begin{array}{l}
i \frac{\partial \zeta}{\partial t} = - \zeta '' +
( |1+\zeta|^2 - 1)(1+\zeta) - u \mu (1+\zeta),\\
\zeta'(t,0)=\zeta'(t,1)=0,\\
\zeta(0,x)=0.
\end{array} \right.
\end{equation}
Proposition \ref{WP-NLS} will be the consequence of the existence and
uniqueness of a weak solution $\zeta$ for (\ref{eq:zeta})
(the conservation of the $L^2$-norm may be proved as in the linear case). 
In order to precise the definition of such a weak solution, let us introduce the 
operator $\mathcal{A}$ defined by
$$\begin{array}{ll}
D(\mathcal{A}):=H^2_{(0)}(0,1),
&
\mathcal{A} \zeta := - \zeta'' + 2 \Re (\zeta).
\end{array}$$
Then for every $\zeta \in H^2_{(0)}(0,1)$ and every $t \in \mathbb{R}$,
we have
$$e^{-i \mathcal{A} t} \zeta =
\sum\limits_{k=0}^{\infty} \left( a_k(t) + i b_k(t) \right) \varphi_k$$
where
$$a_0(t):= \Re (\langle \zeta_0 , \varphi_0 \rangle); \quad
b_0(t):= \Im (\langle \zeta_0 , \varphi_0 \rangle) 
- 2 t \Re (\langle \zeta_0 , \varphi_0 \rangle),$$
$$a_k(t):=
\Re (\langle \zeta_0 , \varphi_k \rangle)
\cos[ \sqrt{\lambda_k(\lambda_k+2)} t ]
+ \sqrt{\frac{\lambda_k}{\lambda_k+2}}
\Im (\langle \zeta_0 , \varphi_k \rangle)
\sin[ \sqrt{\lambda_k(\lambda_k+2)} t ], \forall k \in \mathbb{N}^*,$$
$$b_k(t):=
- \sqrt{\frac{\lambda_k+2}{\lambda_k}}
\Re (\langle \zeta_0 , \varphi_k \rangle)
\sin[ \sqrt{\lambda_k(\lambda_k+2)} t ]
+ \Im (\langle \zeta_0 , \varphi_k \rangle)
\cos[ \sqrt{\lambda_k(\lambda_k+2)} t ], \forall k \in \mathbb{N}^*.$$
Remark that these formulae are only the result of the diagonalization of the matrix $i\mathcal{A}=\left(\begin{array}{cc} 0 &\Delta\\-\Delta +2 & 0 \end{array}\right)$ obtained by the decomposition in real and imaginary part.
Then Proposition \ref{WP-NLS} is equivalent to the following statement.

\begin{Prop} \label{WP-NLS-bis}
Let $\mu \in H^2((0,1),\mathbb{R})$ and $T>0$.
There exists $\delta>0$ such that,
for every $u \in B_\delta [L^2((0,T),\mathbb{R})]$,
there exists a unique weak solution of (\ref{eq:zeta}), i.e.
a function $\zeta \in C^0([0,T],H^2_{(0)})$ such that the
following equality holds in $H^2_{(0)}$ for every $t \in [0,T]$
\begin{equation} \label{WS:zeta}
\zeta(t)=\int_0^t e^{-i \mathcal{A}(t-s)} \Big(
[ |1+\zeta(s)|^2 - 1][1+\zeta(s)] -2 \Re[\zeta(s)] - u(s) \mu [1+\zeta(s)]
\Big) ds.
\end{equation}
\end{Prop}

The proof of Proposition \ref{WP-NLS-bis} relies on the following Lemma.

\begin{Lem} \label{Lem:G-NLS}
Let $T>0$ and $f \in L^2((0,T),H^2)$. The function
$G:t \mapsto \int_0^t e^{-i\mathcal{A}(t-s)} f(s) ds$ belongs to
$C^0([0,T],H^2_{(0)})$, moreover
$$\|G\|_{L^\infty((0,T),H^2_{(0)})} \leqslant c_0(T) \|f\|_{L^2((0,T),H^2)}$$
where the constants $c_0(T)$ are uniformly bounded 
for $T$ lying in bounded intervals.
\end{Lem}

\noindent \textbf{Proof of Lemma \ref{Lem:G-NLS}:} The proof of this Lemma
is similar to the one of Lemma \ref{Lem:G}.  By definition, we have
$$G(t)=\sum\limits_{k=0}^{\infty} \sum\limits_{a=1}^{4} 
\left(\int_0^t y_k^a (t,s) ds\right) \varphi_k$$
where
$$y_k^1(t,s):=
\Re (\langle f(s) , \varphi_k \rangle)
\cos[ \sqrt{\lambda_k(\lambda_k+2)} (t-s) ], \forall k \in \mathbb{N},$$
$$y_k^2(t,s):=
\sqrt{\frac{\lambda_k}{\lambda_k+2}}
\Im (\langle f(s) , \varphi_k \rangle)
\sin[ \sqrt{\lambda_k(\lambda_k+2)} (t-s) ], \forall k \in \mathbb{N}^*,$$
$$y_k^3(t,s):= - i \sqrt{\frac{\lambda_k+2}{\lambda_k}}
\Re (\langle f(s) , \varphi_k \rangle)
\sin[ \sqrt{\lambda_k(\lambda_k+2)} (t-s) ], \forall k \in \mathbb{N}^*,$$
$$y_k^4(t,s):= i \Im (\langle f(s) , \varphi_k \rangle)
\cos[ \sqrt{\lambda_k(\lambda_k+2)} (t-s) ], \forall k \in \mathbb{N},$$
$$y_0^2(t,s):=0, y_0^3(t,s):=-2t\Re (\langle f(s) , \varphi_k \rangle).$$
We have
$$\| G(t) \|_{H^2_{(0)}} \leqslant
\sum\limits_{a=1}^{4} \left(
\sum\limits_{k=1}^{\infty} \Big| k_*^2 \int_0^t y_k^a(t,s) ds \Big|^2
\right)^{1/2}.$$
Let us prove that there exists a constant $c=c(t)>0$ 
(uniformly bounded on bounded intervals of $t$) such that
\begin{equation} \label{ccl-a=1}
\left(
\sum\limits_{k=1}^{\infty} \Big| k_*^2 \int_0^t y_k^1(t,s) ds \Big|^2
\right)^{1/2}
\leqslant
c(t) \| f \|_{L^2((0,t),H^2)}.
\end{equation}
(the other terms may be treated in the same way).
Integrations by part give, for almost every $s \in (0,T)$,
$$\langle f(s) , \varphi_k \rangle
= \frac{\sqrt{2}}{(k\pi)^2} \left(
(-1)^k f'(s,1) - f'(s,0) - \int_0^1 f''(s,x) \cos(k\pi x) dx
\right), \forall k \in \mathbb{N}^*.$$
Thus, we have, for every $k \in \mathbb{N}^*$,
$$\begin{array}{ll}
k^2 \int_0^t y_k^1(t,s) ds=
& \frac{\sqrt{2}(-1)^k}{(\pi)^2} \int_0^t f'(s,1) 
\cos[ \sqrt{\lambda_k(\lambda_k+2)} (t-s) ] ds
\\ &
+ \frac{\sqrt{2}}{(\pi)^2} \int_0^t f'(s,0) 
\cos[ \sqrt{\lambda_k(\lambda_k+2)} (t-s) ] ds
\\ &
- \frac{\sqrt{2}}{(\pi)^2} \int_0^t 
\langle f''(s) , \varphi_k \rangle \cos[ \sqrt{\lambda_k(\lambda_k+2)} (t-s) ] ds.
\end{array}$$
We get (\ref{ccl-a=1}) thanks to Corollary \ref{Cor3}, as in 
the proof of Lemma \ref{Lem:G}. $\Box$
\\

\noindent \textbf{Proof of Proposition \ref{WP-NLS-bis}:}
We introduce the function
$g:\mathbb{C} \rightarrow \mathbb{C}$ defined by $g(z):=[|1+z|^2-1][1+z]$. 
We have $dg(0).\zeta =2 \Re (\zeta)$. 
Let $c_0=c_0(T)$ be as in Lemma \ref{Lem:G-NLS}. Let $c_1, c_2, c_3 >0$ be such that
\begin{equation} \label{def:c1}
\| g(\zeta) - dg(0).\zeta \|_{H^2}
\leqslant c_1 \Big( \|\zeta\|_{H^2_{(0)}}^2 + \|\zeta\|_{H^2_{(0)}}^3 \Big),
\forall \zeta \in H^2_{(0)},
\end{equation}
\begin{equation} \label{def:c2}
\| g(\tilde{\zeta}) - g(\zeta) - dg(0).(\tilde{\zeta} - \zeta) \|_{H^2} 
\leqslant c_2 \|\zeta-\tilde{\zeta}\|_{H^2_{(0)}}
\max\{ \| \xi \|_{H^2_{(0)}} , \| \xi \|_{H^2_{(0)}}^2 ;
\xi \in \{ \zeta , \tilde{\zeta} \} \}, 
\forall \zeta, \tilde{\zeta} \in H^2_{(0)},
\end{equation}
\begin{equation} \label{def:c3}
\| \mu \zeta \|_{H^2} \leqslant c_3 \|\zeta\|_{H^2_{(0)}},
\forall \zeta \in H^2_{(0)}.
\end{equation}
Let $R>0$ be small enough so that
\begin{equation} \label{def:R}
c_0 c_1 \sqrt{T} (R^2 + R^3) < \frac{R}{2}
\text{  and  }
c_0 c_2 \sqrt{T} \max\{ R , R^2 \} < \frac{1}{4}.
\end{equation}
Let $\delta>0$ be small enough so that
\begin{equation} \label{def:delta}
c_0 \delta c_3  (1+R) < \frac{R}{2}
\text{ and }
c_0 \delta c_3 < \frac{1}{4}.
\end{equation}
Let $u \in L^2((0,T),\mathbb{R})$ be such that $\|u\|_{L^2(0,T)} < \delta$.
We consider the map
$$\begin{array}{cccc}
F: & \overline{B}_R[C^0([0,T],H^2_{(0)})] & \rightarrow & 
\overline{B}_R[C^0([0,T],H^2_{(0)})] \\
   &          \zeta             & \mapsto     & \xi
\end{array}$$
where $\xi:=F(\zeta)$ is defined by
$$\xi(t) = -i \int_0^t e^{-i \mathcal{A}(t-s)} \Big(
[ g(\zeta(s)) - dg(0).\zeta(s) - u(s) \mu [1+\zeta(s)]
\Big) ds.$$

For $\zeta \in \overline{B}_R[C^0([0,T],H^2_{(0)})]$, the function
$g(\zeta) - dg(0).\zeta - u \mu [1+\zeta]$
belongs to $L^2((0,T),H^2)$, thus $\xi$ belongs to
$C^0([0,T],H^2_{(0)})$ thanks to Lemma \ref{Lem:G-NLS}. 
Moreover, using (\ref{def:c1}), (\ref{def:c3}), (\ref{def:R}), (\ref{def:delta}), we get
$$\begin{array}{ll}
\| \xi \|_{L^\infty((0,T),H^2_{(0)})}
& \leqslant
c_0 \Big\| g(\zeta) - dg(0).\zeta - u \mu [1+\zeta]\Big\|_{L^2((0,T),H^2)}
\\ & \leqslant
c_0 \Big[
\sqrt{T} \| g(\zeta) - dg(0).\zeta  \|_{L^\infty((0,T),H^2)}
+ \| u\|_{L^2(0,T)} \| \mu [1+\zeta] \|_{L^\infty((0,T),H^2)}
\Big]
\\ & \leqslant
c_0 \Big[
\sqrt{T} c_1 (R^2 + R^3)
+ \delta c_3 (1+R)
\Big]
\\ & \leqslant 
R.
\end{array}$$
Thus, $F$ takes values in $\overline{B}_R[C^0([0,T],H^2_{(0)})]$.

For $\zeta, \tilde{\zeta} \in \overline{B}_R[C^0([0,T],H^2_{(0)})]$, 
using (\ref{def:c2}), (\ref{def:c3}), (\ref{def:R}), (\ref{def:delta}), we get
$$\begin{array}{ll}
\| \xi-\tilde{\xi} \|_{L^\infty((0,T),H^2_{(0)})}
& \leqslant c_0 \Big\|
g(\zeta) - g(\tilde{\zeta}) - dg(0).(\zeta-\tilde{\zeta})
- u \mu (\zeta-\tilde{\zeta})
\Big\|_{L^2((0,T),H^2)}
\\ & \leqslant
c_0 \Big[
\sqrt{T} c_2 \| \zeta-\tilde{\zeta} \|_{L^\infty((0,T),H^2_{(0)})} \max\{R,R^2\} 
+ \delta c_3 \| \zeta - \tilde{\zeta}\|_{L^\infty((0,T),H^2_{(0)})}
\Big]
\\ & \leqslant 
\frac{1}{2}  \| \zeta - \tilde{\zeta}\|_{L^\infty((0,T),H^2_{(0)})}.
\end{array}$$
Thus $F$ is a contraction. $\Box$

\subsection{$C^1$-regularity of the end-point map}

Let $T>0$ and $\delta>0$ be as in Proposition \ref{WP-NLS}. 
Let 
$$V_T:= \left\{ \varphi \in L^2(0,1) ; \Re \left( e^{iT} \int_0^1 \varphi(x) dx \right)=0
\right\},$$
and $P_T:L^2(0,1) \rightarrow V_T$ be the associated orthogonal projection.
Then, the following map is well defined
\begin{equation} \label{def:Theta-NLS}
\begin{array}{cccc}
\Theta_T: & B_{\delta}[L^2((0,T),\mathbb{R})] & \rightarrow & H^2_{(0)}(0,1) \\
          &          u           & \mapsto     & P_T[\psi(T)],
\end{array}
\end{equation}
where $\psi$ solves (\ref{NLS}), (\ref{IC-NLS}).
We want to prove that the map $\Theta_T$ is $C^1$ on a neighborhood of zero.
We have seen that $\psi(t)=e^{-it}(1+\zeta(t))$, where $\zeta$ solves
(\ref{eq:zeta}). Thus, it is sufficient to prove the following statement.

\begin{Prop} \label{Prop:C1-NLS}
Let $\mu \in H^2((0,1),\mathbb{R})$, $T>0$, 
$\delta$ be as in Proposition \ref{WP-NLS-bis},
and
$$\begin{array}{cccc}
\tilde{\Theta}_T: & B_\delta [L^2((0,T),\mathbb{R})] & \rightarrow &  H^2_{(0)}(0,1)
\\
& u & \mapsto & \zeta(T) ,
\end{array}$$
where $\zeta$ solves (\ref{eq:zeta}).
There exists $\delta' \in (0,\delta)$ such that the map $\Theta_T$ is
$C^1$ on $B_{\delta'}[L^2((0,T),\mathbb{R})]$. Moreover, 
for every $u \in B_{\delta'}[L^2((0,T),\mathbb{R})]$ and
$v \in L^2((0,T),\mathbb{R})$
we have 
\begin{equation} \label{diff-ThetaT-NLS}
d\tilde{\Theta}_T(u).v=\xi(T)
\end{equation} 
where $\xi$ solves
\begin{equation} \label{eq:xi}
\left\lbrace \begin{array}{l}
i \frac{\partial \xi}{\partial t} =
- \xi'' + dg(\zeta).\xi - u \mu \xi - v \mu (1+\zeta), \\
\xi'(t,0)=\xi'(t,1)=0,\\
\xi(0,x)=0,
\end{array}\right.
\end{equation}
and $\zeta$ solves (\ref{eq:zeta}).
\end{Prop}

\noindent \textbf{Proof of Proposition \ref{Prop:C1-NLS}:}
We use the same notations $c_0, c_1, c_2, c_3, R, \delta$ as in the proof
of Proposition \ref{WP-NLS-bis}, in particular, the relations
(\ref{def:c1}), (\ref{def:c2}), (\ref{def:c3}), (\ref{def:R}), (\ref{def:delta})
are satisfied. We introduce constants $c_4, c_5>0$ such that
\begin{equation} \label{def:c4}
\| [ dg(\zeta) - dg(0) ].h \|_{H^2} \leqslant c_4 \|h\|_{H^2_{(0)}}
\max\{ \|\zeta\|_{H^2_{(0)}}, \|\zeta\|_{H^2_{(0)}}^2 \},
\forall \zeta, h \in H^2_{(0)}, 
\end{equation}
\begin{equation} \label{def:c5}
\| g(\tilde{\zeta}) - g(\zeta) - dg(0).(\tilde{\zeta} - \zeta ) \|_{H^2}
\leqslant c_5 \| \tilde{\zeta} - \zeta \|_{H^2_{(0)}} 
\max\{ \| \xi \|_{H^2_{(0)}} , \| \xi \|_{H^2_{(0)}}^2 ;
\xi \in \{ \zeta , \tilde{\zeta} \} \},
\forall \zeta, \tilde{\zeta} \in H^2_{(0)}.
\end{equation}
Moreover, we assume that
\begin{equation} \label{def:R-bis}
c_0 \sqrt{T} \max \{ c_4 , c_5 \} \max \{R,R^2 \} < \frac{1}{4}
\end{equation}
(this additional assumption may change $\delta$ into a smaller value $\delta'$).
\\

Let $u,v \in B_\delta[L^2((0,T),\mathbb{R})]$ be such that $(u+v) \in B_\delta[L^2(0,T)]$.
Let $\zeta$, $\xi$ and $\tilde{\zeta}$ be the solutions of
(\ref{eq:zeta}), (\ref{eq:xi}) and
$$\left\lbrace \begin{array}{l}
i \frac{\partial \tilde{\zeta}}{\partial t} = - \tilde{\zeta} '' +
( |1+\tilde{\zeta}|^2 - 1)(1+\tilde{\zeta}) - (u+v) \mu (1+\tilde{\zeta}),\\
\tilde{\zeta}'(t,0)=\tilde{\zeta}'(t,1)=0,\\
\tilde{\zeta}(0,x)=0.
\end{array} \right.$$
The existence of $\xi$ may be proved in a similar way as the existence of $\zeta$.
\\

\emph{First step: Let us prove that}
\begin{equation} \label{cont-NLS}
\| \tilde{\zeta} - \zeta  \|_{L^\infty((0,T),H^2_{(0)})}
\leqslant 2 c_0 c_3 \|1+\zeta\|_{L^\infty((0,T),H^2_{(0)})} \|v\|_{L^2}.
\end{equation}
Thanks to Lemma \ref{Lem:G-NLS}, (\ref{def:c5}),(\ref{def:c3}), (\ref{def:R-bis})
and (\ref{def:delta}), we have
$$\begin{array}{ll}
&
\| \zeta-\tilde{\zeta} \|_{L^\infty((0,T),H^2_{(0)})}
\\ 
\leqslant & 
c_0
\Big\| 
g(\tilde{\zeta}) - g(\zeta) - dg(0).(\tilde{\zeta}-\zeta)
- (u+v) \mu (\tilde{\zeta}-\zeta) - v \mu (1+\zeta)
\Big\|_{L^2((0,T),H^2)}
\\  \leqslant &
c_0 \Big[
\sqrt{T} c_5 \| \tilde{\zeta} - \zeta \|_{L^\infty((0,T),H^2_{(0)})} \max \{ R, R^2 \} 
+ \delta c_3 \| \tilde{\zeta} - \zeta \|_{L^\infty((0,T),H^2_{(0)})}
\\ &
+ \|v\|_{L^2} c_3 \|1+\zeta\|_{L^\infty((0,T),H^2_{(0)})}
\Big]
\\ \leqslant &
\frac{1}{2} \| \zeta-\tilde{\zeta} \|_{L^\infty((0,T),H^2_{(0)})}
+ c_0  \|v\|_{L^2} c_3 \|1+\zeta\|_{L^\infty((0,T),H^2_{(0)})}
\end{array}$$
which gives (\ref{cont-NLS}).
\\

\emph{Second step: Let us prove that the linear map
$$\begin{array}{|ccc}
L^2(0,T) & \rightarrow & H^2_{(0)}(0,1) \\
 v       & \mapsto     & \xi(T)
\end{array}$$
is continuous.} Thanks to Lemma \ref{Lem:G-NLS}, (\ref{def:c4}), (\ref{def:c3}),
(\ref{def:R-bis}) and (\ref{def:delta}), we have
$$\begin{array}{ll}
\| \xi \|_{L^\infty((0,T),H^2_{(0)})}
& \leqslant c_0 
\Big\|
[dg(\zeta)-dg(0)].\xi - u \mu \xi - v \mu (1+\zeta)
\Big\|_{L^2((0,T),H^2)}
\\ & \leqslant c_0 \Big[
\sqrt{T} c_4 \|\xi\|_{L^\infty((0,T),H^2_{(0)})} \max \{R,R^2\} +
\delta c_3 \|\xi\|_{L^\infty((0,T),H^2_{(0)})} 
\\ & +
\|v\|_{L^2} c_3 \|1+\zeta\|_{L^\infty((0,T),H^2_{(0)})}
\Big]
\\ & \leqslant
\frac{1}{2} \| \xi \|_{L^\infty((0,T),H^2_{(0)})}
+ c_0 \|v\|_{L^2} c_3 \|1+\zeta\|_{L^\infty((0,T),H^2_{(0)})},
\end{array}$$
which gives
\begin{equation} \label{majo:xi-NLS}
\|\xi\|_{L^\infty((0,T),H^2_{(0)})} \leqslant
2 c_0  c_3  \|v\|_{L^2} \|1+\zeta\|_{L^\infty((0,T),H^2_{(0)})}.
\end{equation}

\emph{Third step: Let us prove that $\tilde{\Theta}_T$ is differentiable 
and that (\ref{diff-ThetaT-NLS}) holds.}
Let $\Delta:=\tilde{\zeta}-\zeta-\xi$. We want to prove that
$$\| \Delta(T) \|_{H^2_{(0)}} = o ( \|v\|_{L^2} )
\text{ when } \|v\|_{L^2} \rightarrow 0.$$
Let $\epsilon>0$. There exists $\eta>0$ such that,
for every $f \in L^\infty((0,T),H^2_{(0)})$ with
$\| f-\zeta \|_{L^\infty((0,T),H^2_{(0)})} < \eta$, we have
$$\| g(f) - g(\zeta) - dg(\zeta).(f-\zeta) \|_{L^\infty((0,T),H^2_{(0)})}
< \epsilon \| f - \zeta \|_{L^\infty((0,T),H^2_{(0)})}.$$
Let us assume that $v$ is small enough so that
$$2 c_0 c_3 \|1+\zeta\|_{L^\infty((0,T),H^2_{(0)})} \|v\|_{L^2}
< \eta.$$
Then, thanks to Lemma \ref{Lem:G-NLS} and (\ref{cont-NLS}), 
(\ref{def:c4}) and (\ref{def:c3}), we have
$$\begin{array}{ll}
& \| \Delta \|_{L^\infty((0,T),H^2_{(0)})}
\\ \leqslant &
c_0 \Big\|
g(\tilde{\zeta}) - g(\zeta) - dg(\zeta).(\tilde{\zeta}-\zeta)
+ [dg(\zeta)-dg(0)].\Delta - (u+v)\mu \Delta - v \mu \xi
\Big\|_{L^2((0,T),H^2)}
\\  \leqslant &
c_0 \Big[
\sqrt{T} \epsilon \|\tilde{\zeta}-\zeta\|_{L^\infty((0,T),H^2_{(0)})} 
+
\sqrt{T} c_4 (R+R^2) \| \Delta \|_{L^\infty((0,T),H^2_{(0)})}
\\ & 
+ \delta c_3 \| \Delta \|_{L^\infty((0,T),H^2_{(0)})}
+
\|v\|_{L^2} c_3 \| \xi \|_{L^\infty((0,T),H^2_{(0)})}
\Big].
\end{array}$$
Thanks to (\ref{def:R-bis}) and (\ref{def:c3}), we get
$$\| \Delta \|_{L^\infty((0,T),H^2_{(0)})}
\leqslant 2 c_0 \Big[
\sqrt{T} \epsilon \|\tilde{\zeta}-\zeta\|_{L^\infty((0,T),H^2_{(0)})} 
+
\|v\|_{L^2} c_3 \| \xi \|_{L^\infty((0,T),H^2_{(0)})}
\Big],$$
which gives the conclusion, thanks to (\ref{cont-NLS}) and (\ref{majo:xi-NLS}). 
\\

The continuity of the map $d\tilde{\Theta}_T$ may be proved with similar arguments. $\Box$

\subsection{Controllability of the linearized system}

The goal of this section is the proof of the following result.

\begin{Prop} \label{Cont-Lin-NLS}
Let $T>0$ and $\mu \in H^2((0,1),\mathbb{R})$ be such that (\ref{hyp_mu-NLS}) holds.
Let $\delta>0$ be as in Proposition \ref{WP-NLS} and
$\Theta_T$ be defined by (\ref{def:Theta-NLS}).
The linear map $d\Theta_T(0):L^2((0,T),\mathbb{R}) \rightarrow V_T \cap H^2_{(0)}(0,1)$
has a continuous right inverse
$d\Theta_T(0)^{-1}:V_T \cap H^2_{(0)}(0,1) \rightarrow L^2((0,T),\mathbb{R})$.
\end{Prop}

\noindent \textbf{Proof of Proposition \ref{Cont-Lin-NLS}:}
It is equivalent to prove that the continuous linear map
$d\tilde{\Theta}_T(0):L^2((0,T),\mathbb{R}) \rightarrow \tilde{V} \cap H^2_{(0)}(0,1)$
has a continuous right inverse, where
$$\tilde{V}:=\left\{ \varphi \in L^2(0,1) ; \Re \int_0^1 \varphi(x) dx =0 \right\}.$$
We have $d\tilde{\Theta}_T(0).v=\xi(T)$ where $\xi$ is the weak solution of
$$\left\lbrace\begin{array}{l}
i \frac{\partial \xi}{\partial t} = - \xi'' 
+ 2 \Re ( \xi )  - v(t) \mu(x) , x \in (0,1), t \in (0,T),
\\
\xi'(t,0)=\xi'(t,1)=0,
\\
\xi(0,x)=0.
\end{array}\right.$$
In particular, we have
$$\xi(T)=i \int_0^T e^{-i \mathcal{A}(T-s)} v(s) \mu ds
= i \sum\limits_{k=0}^{\infty} [a_k(T)+ib_k(T)]\varphi_k$$
where
$$a_0(T)=\langle \mu , \varphi_0 \rangle \int_0^T v(s) ds,$$
$$b_0(T)=-2 \langle \mu , \varphi_0 \rangle \int_0^T (T-t) v(s) ds,$$
$$a_k(T)=\langle \mu , \varphi_k \rangle 
\int_0^T v(s) \cos[ \sqrt{\lambda_k(\lambda_k+2)}(T-s) ] ds,
\forall k \in \mathbb{N}^*,$$
$$b_k(T)=- \sqrt{\frac{\lambda_k+2}{\lambda_k}} \langle \mu , \varphi_k \rangle
\int_0^T v(s) \sin[ \sqrt{\lambda_k(\lambda_k+2)}(T-s) ] ds,
\forall k \in \mathbb{N}^*.$$
For $\xi_f \in \tilde{V} \cap H^2_{(0)}(0,1)$,
the equality $\xi(T)=\xi_f$ is equivalent to the following 
trigonometric moment problem 
$$\left\lbrace \begin{array}{l}
\int_0^T v(s) ds = d_0(\xi_f):=\Im
\frac{\langle \xi_f , \varphi_0 \rangle}{\langle \mu , \varphi_0 \rangle},
\\
\int_0^T v(s) e^{-i \sqrt{\lambda_k(\lambda_k+2)} s} ds
=d_k(\xi_f):=
\frac{e^{-i \sqrt{\lambda_k(\lambda_k+2)} T}}{\langle \mu , \varphi_k \rangle}
\Big(
\Im \langle \xi_f , \varphi_k \rangle
+ i \sqrt{\frac{\lambda_k}{\lambda_k+2}} \Re  \langle \xi_f , \varphi_k \rangle
\Big), \forall k \in \mathbb{N}^*,
\\
\int_0^T s v(s) ds = \tilde{d}(\xi_f):=T d_0(\xi_f).
\end{array}\right.$$
We conclude thanks to Corollary \ref{Cor:haraux} (in Appendix B). $\Box$
\\

The proof of Theorem \ref{thm_controlNLS} is completed using the same 
arguments as in Section \ref{subsect_proofthm_control_lin} using the inverse 
mapping theorem and the conservation of the $L^2$ norm.

\begin{rk} \label{focusing}
With the same method, one may prove the local exact controllability of the
focusing nonlinear Schr\"odinger equation
$$\left\lbrace \begin{array}{l}
i \frac{\partial \psi}{\partial t}(t,x) = - \frac{\partial^2 \psi}{\partial x^2}(t,x)
- |\psi|^2 \psi(t,x) - u(t) \mu(x) \psi(t,x), x \in (0,1), t \in (0,T),\\
\frac{\partial \psi}{\partial x}(t,0)=\frac{\partial \psi}{\partial x}(t,1)=0,
\end{array}\right.$$
around the reference trajectory
$(\psi_{ref}(t,x)=e^{it},u_{ref}(t)=0)$.
The only difference in the proof is that we get the frequencies
$\sqrt{\lambda_k(\lambda_k-2)}$ (instead of $\sqrt{\lambda_k(\lambda_k+2)}$)
in the moment problem. When the space domain is the interval $(0,1)$, 
then all the quantities $\lambda_k(\lambda_k-2)$, for $k \in \mathbb{N}^*$, are positive
(because $\lambda_k=(k\pi)^2$), thus there is no additional difficulty.
When the space domain is different, for instance $(0,a)$ with $a$ large,
then $\lambda_{k}=(k\pi/a)^2$, thus a finite number of the quantities
$\lambda_k(\lambda_k-2)$ are negative: we get a new moment problem with
a finite number of moments with real valued exponentials,
and a infinite number of trigonometric moments,
that can be easily solved by adapting the tools used in this article.
\end{rk}

\section{Nonlinear wave equations}
\label{sectNLW}

In this section, we study the nonlinear wave equation 
with Neumann boundary conditions (\ref{wave}).
The goal is the proof of Theorem \ref{thm:ondes}.
In all this section, we use the notations defined in
(\ref{def:A-NLS}), (\ref{vep-NLS}), (\ref{def:Hs-NLS}), (\ref{def:k*})
and all the functions are real valued.
\\

First, let us check that the Cauchy problem is well posed
in $H^3_{(0)} \times H^2_{(0)}(0,1)$, when $u \in L^2(0,T)$.
In order to write the system (\ref{wave}) in first order form,
let us introduce 
\begin{equation} \label{def:AA}
\begin{array}{ll}
D(\mathcal{A}):=H^2_{(0)} \times H^1(0,1),
&
\mathcal{A}:=\left( \begin{array}{cc}
0 & Id \\ -A & 0
\end{array} \right),
\\
D(\mathcal{B}):=L^2 \times L^2 (0,1),
&
\mathcal{B}:=\mu(x) \left( \begin{array}{cc}
0 & 0 \\ Id & Id
\end{array} \right)
\end{array}
\end{equation}
and $F:\mathbb{R}^2 \rightarrow \mathbb{R}^2$ defined by $F(y_1,y_2):=(0,f(y_1,y_2))$.
The operator $\mathcal{A}$ generates a $C^0$-group of bounded operators of
$H^2_{(0)} \times H^1(0,1)$ defined by
$$e^{\mathcal{A} t} \left( \begin{array}{c}
w_0 \\ \dot{w}_0 
\end{array} \right) =
\left( \begin{array}{l}
w(t) \\ \dot{w}(t)
\end{array} \right),$$
where
$$w(t)=  \left( \langle w_0 , \varphi_0 \rangle + \langle \dot{w}_0 , \varphi_0 \rangle t
\right) \varphi_0 
+
\sum\limits_{k=1}^{\infty} 
\left( \langle w_0 , \varphi_k \rangle \cos(\sqrt{\lambda_k} t) +
\frac{1}{\sqrt{\lambda_k}} \langle \dot{w}_0 , \varphi_k \rangle \sin(\sqrt{\lambda_k}t)
\right) \varphi_k,$$
$$\dot{w}(t)=\langle \dot{w}_0 , \varphi_0 \rangle \varphi_0 
+
\sum\limits_{k=1}^{\infty} 
\left( - \sqrt{\lambda_k} \langle w_0 , \varphi_k \rangle \sin(\sqrt{\lambda_k} t) +
\langle \dot{w}_0 , \varphi_k \rangle \cos(\sqrt{\lambda_k}t)
\right) \varphi_k.$$
With the notation
$$\begin{array}{ll}
\mathcal{W}:=\left(  \begin{array}{c}
w \\ \frac{\partial w}{\partial t}
\end{array} \right),
&
\mathcal{W}_0:=\left(  \begin{array}{c}
1 \\ 0
\end{array} \right),
\end{array}$$
the equation (\ref{wave}) may be written
\begin{equation} \label{wave-1st}
\frac{\partial \mathcal{W}}{\partial t}(t,x) =
\mathcal{A} \mathcal{W}(t,x) + 
F(\mathcal{W})+ u(t) \mathcal{B} \mathcal{W}(t,x), x \in (0,1).
\end{equation}

\begin{Prop} \label{WP-wave}
Let $\mu \in H^2(0,1)$, $T>0$, $f \in C^3(\mathbb{R}^2,\mathbb{R})$ be such that
$f(1,0)=0$ and $\nabla f(1,0)=0$.
There exists $\delta>0$ such that, for every $u \in B_\delta[L^2(0,T)]$,
there exists a unique weak solution of (\ref{wave-1st}), (\ref{IC-wave}), i.e. a function
$\mathcal{W} \in C^0([0,T],H^3_{(0)} \times H^2_{(0)})$ such that 
the following equality holds in $H^3_{(0)} \times H^2_{(0)}(0,1)$,
for every $t \in [0,T]$,
\begin{equation} \label{SolutionFaibleonde}
\mathcal{W}(t) = e^{\mathcal{A}t} \mathcal{W}_0
+ \int_0^t e^{\mathcal{A}(t-\tau)} 
\Big( F(\mathcal{W}(\tau)) + u(\tau) \mathcal{B} \mathcal{W}(\tau) 
+ \mathcal{F}(\tau) \Big) d\tau.
\end{equation}
\end{Prop}

The proof of this proposition relies on the following Lemma.

\begin{Prop} \label{Prop-wave}
Let $T>0$ and $g \in L^2((0,T),H^2)$. The function $G$ defined by
$$G(t):=\int_0^t e^{\mathcal{A} s} 
\left( \begin{array}{c}
0 \\ g(s)
\end{array} \right) ds$$
belongs to $C^0([0,T],H^3_{(0)} \times H^2_{(0)})$.
Moreover, there exists a constant $c_0(T)>0$,
uniformly bounded for $T$ lying in bounded intervals,
such that, for every $g \in L^2((0,T),H^2)$,
\begin{equation} \label{majo_lem_wave}
\| G \|_{L^\infty((0,T),H^3_{(0)} \times H^2_{(0)})}
\leqslant c_0(T) \|g\|_{L^2((0,T),H^2)}.
\end{equation}
\end{Prop}

\noindent \textbf{Proof of Proposition \ref{Prop-wave}:}
We have, for every $t \in [0,T]$,
$$G(t)=\int_0^t
\left(
\begin{array}{l}
\langle g(s) , \varphi_0 \rangle s \varphi_0 
+\sum\limits_{k=1}^{\infty} 
\frac{\langle g(s) , \varphi_k \rangle}{\sqrt{\lambda_k}}  \sin(\sqrt{\lambda_k}s) \varphi_k
\\
\langle g(s) , \varphi_0 \rangle \varphi_0 
+
\sum\limits_{k=1}^{\infty} 
\langle g(s) , \varphi_k \rangle \cos(\sqrt{\lambda_k}s) \varphi_k
\end{array}
\right) ds.$$
Thus, there exists $C>0$ such that
$$\|G(t)\|_{H^3_{(0)} \times H^2_{(0)}}
\leqslant C 
\sum\limits_{k=0}^{\infty} 
\Big| k_*^2
\int_0^t \langle g(s) , \varphi_k \rangle e^{i k \pi s} ds
\Big|^2.$$
We get the conclusion as in the previous sections. $\Box$
\\

\noindent \textbf{Proof of Proposition \ref{WP-wave}:}
Let us introduce the constants $c_1, c_2, c_3$ such that
\begin{equation} \label{def:c1-wave}
\| f(w,w_t) \|_{H^2} \leqslant
c_1 \| (w-1,w_t) \|_{H^3_{(0)} \times H^2_{(0)}}^2,
\forall (w,w_t) \in (1,0)+ B_1[H^3_{(0)} \times H^2_{(0)}],
\end{equation}
\begin{equation} \label{def:c2-wave}
\begin{array}{ll}
& \| f(w,w_t) - f(\tilde{w},\tilde{w}_t) \|_{H^2} 
\\ \leqslant &
c_2 \| (w-\tilde{w},w_t-\tilde{w}_t) \|_{H^3_{(0)} \times H^2_{(0)}}
\max\{ \| (w-1,w_t) \|_{H^3_{(0)} \times H^2_{(0)}},
\| (\tilde{w}-1,\tilde{w}_t) \|_{H^3_{(0)} \times H^2_{(0)}}  \},
\\ &
\forall (w,w_t), (\tilde{w},\tilde{w}_t) \in (1,0)+ B_1[H^3_{(0)} \times H^2_{(0)}],
\end{array}
\end{equation}
and (\ref{def:c3}) holds. Let $R \in (0,1)$ be small enough so that
\begin{equation} \label{def:R-wave}
\begin{array}{ll}
\sqrt{T} c_0 c_1 R^2 \leqslant \frac{R}{2},
&
\sqrt{T} c_0 c_2 R \leqslant \frac{1}{4} 
\end{array}
\end{equation}
Let $\delta>0$ be small enough so that 
\begin{equation} \label{def:delta-wave}
\delta c_0 c_3  < \frac{1}{4}, \quad \delta c_0 c_3 (1+R) < \frac{R}{2}.
\end{equation}
Let $u \in B_\delta[L^2(0,T)]$. We consider the map
$$\begin{array}{crcl}
\mathcal{F} : & (1,0)+ \overline{B}_R[C^0([0,T],H^3_{(0)} \times H^2_{(0)})] & 
\rightarrow & (1,0) + \overline{B}_R[C^0([0,T],H^3_{(0)} \times H^2_{(0)})]
\\
& \zeta & \mapsto & \xi
\end{array}$$
where
$$\xi(t)=e^{\mathcal{A}t} \mathcal{W}_0
+ \int_0^t  e^{\mathcal{A}(t-\tau)} 
\Big( F(\zeta(\tau)) + u(\tau) \mathcal{B} \zeta(\tau)  \Big) d\tau,
\forall t \in [0,T].$$

For $\zeta=(w,w_t) \in (1,0) + \overline{B}_R [ C^0([0,T],H^3_{(0)} \times H^2_{(0)}) ]$,
the second component of
$F(\zeta) + u \mathcal{B} \zeta$
belongs to $L^2((0,T),H^2)$, thus $\xi$ belongs to
$C^0([0,T],H^3_{(0)} \times H^2_{(0)})$ thanks to Proposition \ref{Prop-wave}.
Moreover, thanks to (\ref{majo_lem_wave}), (\ref{def:c1-wave}), (\ref{def:c3}), 
(\ref{def:R-wave}), and (\ref{def:delta-wave}), we have, for every $t \in [0,T]$,
$$\begin{array}{ll}
\| \xi(t)-(1,0) \|_{H^3_{(0)} \times H^2_{(0)}}
& \leqslant 
c_0 \| f(w,w_t) + u \mu [w+w_t] \|_{L^2((0,T),H^2)}
\\ & \leqslant
c_0 \Big(
\sqrt{T} \|f(w,w_t)\|_{L^\infty((0,T),H^2)}
+ \|u\|_{L^2(0,T)} \| \mu [w+w_t] \|_{L^\infty((0,T),H^2)}
\Big)
\\ & \leqslant
c_0 ( \sqrt{T} c_1 R^2 + \delta c_3 (R+1) )
\\ & \leqslant R.
\end{array}$$
Thus, $\mathcal{F}$ takes values in 
$(1,0)+\overline{B}_R[C^0([0,T],H^3_{(0)} \times H^2_{(0)})]$.

For $\zeta=(w,w_t), \tilde{\zeta}=(\tilde{w},\tilde{w}_t) 
\in (1,0)+\overline{B}_R[C^0([0,T],H^3_{(0)} \times H^2_{(0)})]$,
thanks to (\ref{def:c2-wave}), (\ref{def:c3}), (\ref{def:R-wave}) and
(\ref{def:delta-wave}), we have
$$\begin{array}{ll}
& \| \mathcal{F}(\zeta) - \mathcal{F}(\tilde{\zeta}) 
\|_{L^\infty((0,T),H^3_{(0)} \times H^2_{(0)})}
\\ \leqslant & 
c_0 \| f(w,w_t) - f(\tilde{w},\tilde{w}_t) + u \mu [w-\tilde{w}+w_t-\tilde{w}_t] 
\|_{L^2((0,T),H^2)}
\\  \leqslant &
c_0 \Big[ \sqrt{T} \| f(w,w_t) - f(\tilde{w},\tilde{w}_t) \|_{L^\infty((0,T),H^2)}
+ \delta c_3 \| \zeta - \tilde{\zeta} \|_{L^\infty((0,T),H^3_{(0)} \times H^2_{(0)})}
\Big]
\\  \leqslant & 
c_0 \Big[ \sqrt{T} c_2 R 
\| \zeta - \tilde{\zeta} \|_{L^\infty((0,T),H^3_{(0)} \times H^2_{(0)})}
+ \delta c_3 \| \zeta - \tilde{\zeta} \|_{L^\infty((0,T),H^3_{(0)} \times H^2_{(0)})}
\Big]
\\ \leqslant &
\frac{1}{2} \| \zeta - \tilde{\zeta} \|_{L^\infty((0,T),H^2_{(0)})}
\end{array}$$
Thus $\mathcal{F}$ is a contraction. $\Box$
\\

Let $T>0$, $\mu \in H^2(0,1)$,
$f \in C^3(\mathbb{R}^2,\mathbb{R})$ be such that $f(1,0)=0$,  $\nabla f(1,0)=0$
and $\delta>0$ be as in Proposition \ref{WP-wave}.
Then, the following map is well defined
\begin{equation} \label{def:ThetaT-waves}
\begin{array}{cccc}
\Theta_T: & B_\delta[L^2(0,T)] & \rightarrow & H^3_{(0)} \times H^2_{(0)} \\
          &    u               & \mapsto     & (w,w_t)(T)
\end{array}
\end{equation}
where $(w,w_t)$ is the weak solution of (\ref{wave}), (\ref{IC-wave}).
Working as in the previous section, one may prove the following statements.

\begin{Prop} \label{C1-wave}
Let $\mu \in H^2(0,1)$, $T>0$, 
$f \in C^3(\mathbb{R}^2,\mathbb{R})$ be such that $f(1,0)=0$,  $\nabla f(1,0)=0$
and $\delta>0$ be as in Proposition  \ref{WP-wave}.
The map $\Theta_T$ defined by (\ref{def:ThetaT-waves}) is $C^1$.
Moreover, for every $u \in B_\delta[L^2(0,T)]$ and
$v \in L^2(0,T)$, we have $d\Theta_T(u).v=(W,W_t)(T)$, where
$(W,W_t)$ is the weak solution of
\begin{equation} \label{lin-waves}
\left\lbrace \begin{array}{l}
W_{tt}=W_{xx}+\frac{\partial f}{\partial y_1}(w,w_t).W
+ \frac{\partial f}{\partial y_2}(w,w_t).W_t + u(t) \mu [W+W_t]
+ v(t) \mu(x) [ w+w_t], 
\\
W_x(t,0)=W_x(t,1)=0,
\\
(W,W_t)(0,x)=0,
\end{array} \right.
\end{equation}
and $(w,w_t)$ is the weak solution of (\ref{wave}), (\ref{IC-wave}).
\end{Prop}

\begin{Prop} 
Let $T>2$, $\mu \in H^2(0,1)$ be such that (\ref{hyp_mu-NLS}) holds
and $f \in C^3(\mathbb{R}^2,\mathbb{R})$ be such that $f(1,0)=0$,  $\nabla f(1,0)=0$.
The linear map
$d\Theta_T(0):L^2(0,T) \rightarrow H^3_{(0)} \times H^2_{(0)}$
has a continuous right inverse
$d\Theta_T(0)^{-1}:H^3_{(0)} \times H^2_{(0)} \rightarrow L^2(0,T)$. 
\end{Prop}

The proof is the same except that the gap between the eigenvalues does not tend to infinity and we use Corollary \ref{Corwave}.

\section{Conclusion, open problems, perspectives}
\label{sec:ccl}

In this article, we have proposed a method for the proof of the
local exact controllability for linear and nonlinear
bilinear systems. We have applied it to Schr\"odinger and wave equations,
showing it works for a wide range of problems.
It also works on other equations (for instance it may prove an optimal 
version of the controllability result proved in \cite{KB-poutres}
for a 1D Beam equation).
\\

In this article, we have presented various examples of application of the method. 
However, they all have in common that the linearized system fulfills a gap condition 
on the eigenvalues of the operator. This condition is not necessarily realized 
for the Schr\"odinger equation in higher space dimensions. 
Even in two dimension, we do not know any example of domain where 
it is true. So, one challenging question is the extension (or the impossibility 
to do it) of these results to other dimensions.

\appendix

\section{Genericity of the assumption on $\mu$}
\label{appendix:hyp_mu}

The goal of this section is the proof of the following result.

\begin{Prop}
The set
$\{ \mu \in H^{3}((0,1),\mathbb{R}) ; (\ref{hyp_mu}) \text{ holds} \}$
is dense in $H^{3}((0,1),\mathbb{R})$.
\end{Prop}

\noindent \textbf{Proof:} First, let us notice that
$$\mathcal{V}:=\{ \mu \in H^{3}((0,1),\mathbb{R}) ;
\mu'(1) \pm \mu'(0) \neq 0 \}$$
is a dense open subset of $H^{3}((0,1),\mathbb{R})$.
Now, let us prove that the set
$$\mathcal{U}:=\{ \mu \in \mathcal{V} ;
\langle \mu \varphi_{1} , \varphi_{k} \rangle \neq 0, 
\forall k \in \mathbb{N}^{*} \}$$
is dense in $H^{3}((0,1),\mathbb{R})$.
It is sufficient to prove that this set is dense in $\mathcal{V}$. 
For $n \in \mathbb{N}$, we introduce the set
$$\mathcal{U}_{n} := \{ \mu \in \mathcal{V} ;
\langle \mu \varphi_{1} , \varphi_{k} \rangle \neq 0, 
\forall k \in \{1,...,n\} \},$$
with the convention 
$\mathcal{U}_{0} := \mathcal{V}$.
Then the sequence $(\mathcal{U}_{n})_{n \in \mathbb{N}}$ is decreasing and
$$\mathcal{U}=\bigcap_{n=0}^{\infty} \mathcal{U}_{n}.$$
Thanks to Baire Lemma, it is sufficient to check that, for every $n \in \mathbb{N}$, 
$\mathcal{U}_{n+1}$ is dense in $\mathcal{U}_{n}$
for the $H^{3}((0,1),\mathbb{R})$-topology.
Let $n \in \mathbb{N}$ and let 
$\mu \in \mathcal{U}_{n} - \mathcal{U}_{n+1}$.
Then $\mu \in \mathcal{V}$,
$\langle \mu \varphi_{1} , \varphi_{k} \rangle \neq 0$ for $k=1,...,n$ 
and $\langle \mu \varphi_{1} ,\varphi_{n+1} \rangle = 0$.
Thanks to (\ref{coeff-explicit-x2}), 
$\mu + \epsilon x^{2} \in \mathcal{U}_{n+1}$
for every $\epsilon \in \mathbb{R}$ such that
$$\epsilon \neq - 
\frac{\langle \mu \varphi_{1} , \varphi_{j} \rangle}{\langle x^{2} \varphi_{1} , 
\varphi_{j} \rangle}, \forall j \in \{1,...,n\}.$$
Thus $\mathcal{U}_{n+1}$ is dense in $\mathcal{U}_{n}$.

Finally, thanks to (\ref{hyp_mu-generic}), we have
$$\mathcal{U} \subset 
\{ \mu \in H^{3}((0,1),\mathbb{R}) ; (\ref{hyp_mu}) \text{ holds} \},$$
which gives the conclusion. $\Box$

\begin{Prop} \label{prop_hyp_mu-generic}
The set
$\{ \mu \in H^{3}_{rad}(B^3,\mathbb{R}) ; (\ref{hyp_mu_rad}) \text{ holds} \}$
is dense in $H^{3}(B^3,\mathbb{R})$.
\end{Prop}
\noindent \textbf{Proof:}
We make the same proof. We use the formula 
\begin{eqnarray*}
\left\langle \mu \varphi_1,\varphi_k \right\rangle=\frac{4\pi(-1)^{k+1}}{\lambda_k^{3/2}} \partial_r \mu(s,r=1)-\frac{1}{\lambda_k^2}\int_{B^3} \nabla \Delta (\mu \varphi_1)\cdot \nabla \varphi_k
\end{eqnarray*}
instead of (\ref{hyp_mu-generic}). Moreover, we can find one $\mu(r)=r^2$ that fulfills (\ref{hyp_mu_rad})
$$ \left\langle r^2 \varphi_1,\varphi_k\right\rangle=
2\pi \int_0^1 r^2 sin(\pi r) sin(k\pi r). \Box$$

\section{Moment problems}
\label{App:Moment Pb}

In this section, we recall classical results about moment problems 
(see, for instance \cite{Avdonin}). The proofs are given for sake of completeness.

\subsubsection{Families of vectors in Hilbert spaces}

Let $H$ be a separable Hilbert vector space over $\mathbb{K}= \mathbb{R}$ or $\mathbb{C}$ 
and $\Theta:=(\xi_{j})_{j \in \mathbb{Z}}$ be
a family of vectors of $H$ with $\xi_{j} \neq 0, \forall j \in \mathbb{Z}$.

\begin{Def} \label{Def:minimal}
The family $\Theta$ is minimal in H if, for every $j \in \mathbb{Z}$,
$\xi_{j} \notin \overline{\text{Span} \{ \xi_{i} ; i \in \mathbb{Z}-\{j\} \}}$.
\end{Def}

\begin{Prop} \label{Prop:minimal=biorthogonale}
The family $\Theta$ is minimal in H if and only if
there exists a biorthogonal family $\Theta'=(\xi_{j}')_{j \in \mathbb{Z}}$,
i.e. $\Theta'$ is a family of vectors of $H$ such that
\begin{equation} \label{biorthogonale}
\langle \xi_{i} , \xi_{j}' \rangle = \delta_{i,j}, \forall i,j \in \mathbb{Z}.
\end{equation}
\end{Prop}

\noindent \textbf{Proof of Proposition \ref{Prop:minimal=biorthogonale} : }
We assume $\Theta$ is minimal. For $j \in \mathbb{Z}$, let
$v_{j}$ be the orthogonal projection of $\xi_{j}$ over the
closed vector space $\overline{\text{Span} \{ \xi_{i} , i \neq j \}}$ i.e.
$$v_{j} \in \overline{\text{Span} \{ \xi_{i} , i \neq j \}}
\text{  and  }
\langle \xi_{j} - v_{j} , \xi_{i} \rangle = 0 , \forall i \neq j.$$
Let
$$\xi_{j}':=\frac{\xi_{j}-v_{j}}{\| \xi_{j} - v_{j} \|^{2}} , \forall j \in \mathbb{Z}.$$
Then, the families $(\xi_{j})$ and $(\xi_{j}')$ are biorthogonal.

Now, we assume that there exists a biorthogonal family 
$\Theta'=(\xi_{j}')_{j \in \mathbb{Z}}$.
Let us assume that there exists $j \in \mathbb{Z}$ such that
$\xi_{j} \in \overline{ \text{Span} \{ \xi_{i} ; i \in \mathbb{Z}-\{j\} \} }$.
Then (\ref{biorthogonale}) implies 
$\langle \xi_{j} , \xi_{j}' \rangle =1$
which is a contradiction. $\Box$

\begin{rk}
If $\Theta$ is minimal, then there exists a unique biorthogonal family 
$\Theta'$ such that 
$\Theta' \subset \overline{\text{Span} \{ \xi_{i}  ; i \in \mathbb{Z} \}}$.
In the end of this appendix, the expression
\textquotedblleft the\textquotedblright biorthogonal family of $\Theta$,
refers to this unique biorthogonal family in  
$\overline{\text{Span} \{ \xi_{i}  ; i \in \mathbb{Z} \}}$.
\end{rk}

\begin{Def} \label{Def:Riesz}
The family $\Theta$ is a Riesz basis of $\overline{\text{Span} \Theta}$
if $\Theta$ is the image of some orthonormal family by an isomorphism.
\end{Def}

\begin{rk} \label{Riesz->min}
It is clear that, if $\Theta$ is a Riesz basis of $\overline{\text{Span} \Theta}$, 
then $\Theta$ is minimal in $H$.
\end{rk}

\begin{Prop} \label{Prop:Riesz}

\textbf{(1)} If $\Theta$ is a Riesz basis of $\overline{\text{Span} \Theta}$,
then its biorthogonal family $\Theta'$ 
is also a Riesz basis of $\overline{\text{Span} \Theta}$.

\textbf{(2)} $\Theta$ is a Riesz basis of $\overline{\text{Span} \Theta}$
if and only if there exists $C_{1}, C_{2} \in (0,+\infty)$ such that,
for every scalar sequence $(c_{j})_{j \in \mathbb{Z}}$ with finite support,
\begin{equation} \label{propriete-Riesz}
C_{1} \left( \sum_{j=-\infty}^{\infty} |c_{j}|^{2} \right)^{1/2}
\leqslant
\Big\| \sum_{j=-\infty}^{\infty} c_{j} \xi_{j} \Big\|
\leqslant
C_{1} \left( \sum_{j=-\infty}^{\infty} |c_{j}|^{2} \right)^{1/2}.
\end{equation}

\textbf{(3)} If $\Theta$ is a Riesz basis of $\overline{\text{Span} \Theta}$
then there exists $C>0$ such that, for every $f \in H$, we have
$$\left( \sum\limits_{j \in \mathbb{Z}}
| \langle f , \xi_j \rangle |^2 \right)^{1/2}
\leqslant C \|f\|.$$
\end{Prop}

\noindent \textbf{Proof of Proposition \ref{Prop:Riesz} : }

\textbf{(1)} We assume $\Theta$ is a Riesz basis of $\overline{\text{Span} \Theta}$.
Let $\mathcal{H}$ be an Hilbert space, 
$(\zeta_{j})_{j \in \mathbb{Z}}$ be an orthonormal family of $\mathcal{H}$,
$V : \mathcal{H} \rightarrow \overline{\text{Span} \Theta}$ an isomorphism
such that $\xi_{j} = V(\zeta_{j}), \forall j \in \mathbb{Z}$.
Then the adjoint operator 
$V^{*} : \overline{\text{Span} \Theta} \rightarrow \mathcal{H}$
is also an isomorphism and we have 
$\xi_{j}'=(V^{*})^{-1}(\zeta_{j}), \forall j \in \mathbb{Z}$.
Indeed, for every $j,k \in \mathbb{Z}$, 
$$\delta_{j,k} 
= \langle \xi_{j} , \xi_{k}' \rangle_{H}
= \langle V(\zeta_{j}) , \xi_{k}' \rangle_{H}
= \langle \zeta_{j} , V^{*}(\xi_{k}') \rangle_{\mathcal{H}}.$$
Thus $\Theta'$ is a Riesz basis of $\overline{\text{Span} \Theta}$.

\textbf{(2)} We assume $\Theta$ is a Riesz basis of $\overline{\text{Span} \Theta}$.
Let $\mathcal{H}$ be an Hilbert space, 
$(\zeta_{j})_{j \in \mathbb{Z}}$ be an orthonormal family of $\mathcal{H}$,
$V : \mathcal{H} \rightarrow \overline{\text{Span} \Theta}$ an isomorphism
such that $\xi_{j} = V(\zeta_{j}), \forall j \in \mathbb{Z}$ and
$(c_{j})_{j \in \mathbb{Z}}$ a scalar sequence with finite support.
We have
$$\Big\|         \sum_{j=-\infty}^{\infty} c_{j} \xi_{j}         \Big\|
= \Big\| V \Big[ \sum_{j=-\infty}^{\infty} c_{j} \zeta_{j} \Big] \Big\|
\leqslant \| V \| \Big\| \sum_{j=-\infty}^{\infty} c_{j} \zeta_{j} \Big\|
= \| V \| \left( \sum_{j=-\infty}^{\infty} |c_{j}|^{2} \right)^{1/2}$$
and
$$\left( \sum_{j=-\infty}^{\infty} |c_{j}|^{2} \right)^{1/2}
=\Big\| \sum_{j=-\infty}^{\infty} c_{j} \zeta_{j} \Big\|
=\Big\| V^{-1} \Big[ \sum_{j=-\infty}^{\infty} c_{j} \xi_{j} \Big]  \Big\|
\leqslant \|V^{-1}\| \Big\| \sum_{j=-\infty}^{\infty} c_{j} \xi_{j}   \Big\|,$$
thus, we have (\ref{propriete-Riesz}) with 
$C_{1}=1/\| V^{-1} \|$ and $C_{2}=\|V\|$.

Now, we assume that (\ref{propriete-Riesz}) holds. Then the linear map
$V:l^{2}(\mathbb{Z}, \mathbb{K}) \rightarrow \overline{\text{Span} \Theta}$ 
defined by $V[(c_{j})_{j \in \mathbb{Z}}]=\sum_{j=-\infty}^{\infty} c_{j} \xi_{j}$
is well defined and injective. Let $h \in \overline{\text{Span} \Theta}$.
There exists $(h_{N})_{N \in \mathbb{N}}$ such that $h_{N} \rightarrow h$ 
in $H$ when $N \rightarrow + \infty$ and 
for every $N \in \mathbb{N}$, there exists a sequence 
$c^{(N)}=(c_{j}^{(N)})_{j \in \mathbb{Z}}$ with finite support such that
$h_{N}=\sum_{j=-\infty}^{\infty} c_{j}^{(N)} \xi_{j}$.
Then $(h_{N})_{N \in \mathbb{N}}$ is a Cauchy sequence in $H$, 
thus, thanks to (\ref{propriete-Riesz}), $(c^{(N)})_{N \in \mathbb{N}}$ is a
Cauchy sequence in $l^{2}(\mathbb{Z})$ and there exists $c=(c_{j})_{j \in \mathbb{Z}}
\in l^{2}(\mathbb{N})$ such that $c^{N} \rightarrow c$ in $l^{2}(\mathbb{Z})$.
Then, (\ref{propriete-Riesz}) proves that
$\sum_{j=-\infty}^{\infty} (c_{j}-c_{j}^{(N)}) \xi_{j} \rightarrow 0$ in $H$, i.e.
$h=\sum_{j=-\infty}^{\infty} c_{j} \xi_{j}$. We have proved that $V$ is an isomorphism,
thus  $\Theta$ is a Riesz basis of $\overline{\text{Span} \Theta}$. 

\textbf{(3)} $\overline{\text{Span} \Theta}$ is a close vector subspace
of $H$ thus we have the orthogonal decomposition
$H=\overline{\text{Span} \Theta} + \overline{\text{Span} \Theta}^{\perp}$
and the associated orthogonal projection
$P:H \rightarrow \overline{\text{Span} \Theta}$.
For $f \in H$, we have
$$\begin{array}{ll}
\left( \sum\limits_{j \in \mathbb{Z}}
| \langle f , \xi_j \rangle |^2 \right)^{1/2}
& =
\left( \sum\limits_{j \in \mathbb{Z}}
| \langle Pf , \xi_j \rangle |^2 \right)^{1/2}
\\ & \leqslant
\frac{1}{C_1} \Big\| \sum\limits_{j \in \mathbb{Z}}
\langle Pf , \xi_j \rangle \xi_j' \Big\|
\\ & =
\frac{1}{C_1} \|Pf\|_H \leqslant \frac{1}{C_1} \|f\|. \Box
\end{array}$$

\begin{rk} \label{rk:decomposition}
We have proved that, if $\Theta$ is a Riesz basis of $\overline{\text{Span} \Theta}$,
then, for every $h \in \overline{\text{Span} \Theta}$ there exists 
$c=(c_{j})_{j \in \mathbb{Z}} \in l^{2}(\mathbb{Z},\mathbb{K})$ such that
$h=\sum_{j=-\infty}^{\infty} c_{j} \xi_{j}$. 
Moreover, if $\Theta'$ and $\Theta$ are biorthogonal families, then necessarily
$c_{j}=\langle h , \xi_{j}' \rangle , \forall j \in \mathbb{Z}.$
Thus, every $h \in \overline{\text{Span} \Theta}$ can be decomposed in
the following way
\begin{equation} \label{formule-decomposition}
h=\sum_{j=-\infty}^{\infty} \langle h , \xi_{j}' \rangle \xi_{j}
 =\sum_{j=-\infty}^{\infty} \langle h , \xi_{j} \rangle \xi_{j}'
\end{equation}
where the series converge in $H$ and the coefficients
$( \langle h , \xi_{j}' \rangle )_{j \in \mathbb{Z}}$,
$( \langle h , \xi_{j} \rangle )_{j \in \mathbb{Z}}$,
belong to $l^{2}(\mathbb{Z},\mathbb{K})$. 
\end{rk}

\subsubsection{Abstract moment problems}

Now, we move to the investigation of abstract moment problems:
given a scalar sequence $(d_{j})_{j \in \mathbb{Z}}$ is it possible to find
$f \in H$ such that
$$\langle f , \xi_{j} \rangle = d_{j} , \forall j \in \mathbb{Z}.$$
Let us introduce the operator
$$\begin{array}{cccc}
J_{\Theta} : &  H &  \rightarrow & l^{2}(\mathbb{Z},\mathbb{K}) \\
             &  f & \mapsto      & (\langle f , \xi_{j} \rangle)_{j \in \mathbb{Z}} 
\end{array}$$
with domain
$D_{\Theta} := \{ f \in H ; J_{\Theta}(f) \in l^{2}(\mathbb{Z}) \}$.
It is clear that, if the family $\Theta$ is not complete in $H$, then
the operator $J_{\Theta}$ has a non trivial null space
$\overline{ \text{Span} \Theta }^{\perp}$.
This motivates the introduction of the operator 
$J_{\Theta}^{0} := J_{\Theta} \Big|_{\overline{\text{Span} \Theta}}$.

\begin{Prop} \label{Riesz=solution-moment}
The operator  
$J_{\Theta}^{0} : \overline{\text{Span} \Theta} \rightarrow l^{2}(\mathbb{Z},\mathbb{K})$
is an isomorphism if and only if $\Theta$ is a 
Riesz basis of $\overline{\text{Span} \Theta}$.
\end{Prop}

\noindent \textbf{Proof of Proposition \ref{Riesz=solution-moment} : }
We assume 
$J_{\Theta}^{0} : \overline{\text{Span} \Theta} \rightarrow l^{2}(\mathbb{Z},\mathbb{K})$
is an isomorphism. Let $(\zeta_{j})_{j \in \mathbb{Z}}$ 
be the canonical orthonormal basis of
$l^{2}(\mathbb{Z})$. Then, the family
$$\Big( (J_{\Theta}^{0})^{-1}(\zeta_{j}) \Big)_{j \in \mathbb{Z}}$$
is a Riesz basis of $\overline{\text{Span} \Theta}$. Moreover, it is the
biorthogonal family to $\Theta$ in $\overline{\text{Span} \Theta}$.
Thanks to Proposition \ref{Prop:Riesz} \textbf{(1)}, $\Theta$ is also a  
Riesz basis of $\overline{\text{Span} \Theta}$.

We assume $\Theta$ is a Riesz basis of $\overline{\text{Span} \Theta}$.
Thanks to the Remark \ref{rk:decomposition}, it is clear that
$J_{\Theta}^{0} : \overline{\text{Span} \Theta} \rightarrow l^{2}(\mathbb{Z},\mathbb{K})$
is an isomorphism. $\Box$

\subsubsection{Trigonometric moment problems}

In this section, we recall important results on trigonometric moment problems.
The following Ingham inequality is due to Haraux \cite{Haraux}.

\begin{thm} \label{thm:Haraux}
Let $N \in \mathbb{N}$,
$(\omega_{k})_{k \in \mathbb{Z}}$ be an increasing sequence of real numbers such that
$$\omega_{k+1} - \omega_{k} \geqslant \gamma >0 , 
\forall k \in \mathbb{Z} , |k| \geqslant N,$$
$$\omega_{k+1} - \omega_{k} \geqslant \rho >0 , 
\forall k \in \mathbb{Z},$$
and $T>2\pi / \gamma$.
There exists $C_{1}=C_{1}(\gamma, \rho, N, T), 
C_{2}=C_{2}(\gamma, \rho, N, T) \in (0,+\infty)$ 
such that, for every sequence 
$(c_{k})_{k \in \mathbb{Z}} \in \mathbb{C}^{\mathbb{Z}}$ with finite support, we have
$$C_{1} \sum_{k \in \mathbb{Z}} |c_{k}|^{2}
\leqslant
\int_{0}^{T} \Big| \sum_{k=-\infty}^{+ \infty} c_{k} e^{-i\omega_{k}t} \Big|^{2} dt
\leqslant
C_{2} \sum_{k \in \mathbb{Z}} |c_{k}|^{2}.$$
\end{thm}

Let us introduce the space
$$l^2_r (\mathbb{N},\mathbb{C}):=
\{ (d_k)_{k \in \mathbb{N}} \in l^2(\mathbb{N},\mathbb{C}) ; d_0 \in \mathbb{R} \}.$$
Thanks to Proposition \ref{Prop:Riesz} and Theorem \ref{thm:Haraux}, 
we have the following statement,
which is used in the proof of Proposition \ref{Cont-Lin-H3L2}.

\begin{Cor} \label{Cor:haraux1}
Let $T>0$ and $(\omega_{k})_{k \in \mathbb{N}}$ be an increasing sequence of $[0,+\infty)$
such that $\omega_{0}=0$ and 
$$\omega_{k+1} - \omega_{k} \rightarrow + \infty \text{ when } k \rightarrow + \infty.$$
There exists a continuous linear map
$$\begin{array}{cccc}
L: & l^2_r(\mathbb{N},\mathbb{C}) & \rightarrow & L^2((0,T),\mathbb{R}) \\
   &                 d            & \mapsto     & L(d)
\end{array}$$
such that, for every $d=(d_k)_{k \in \mathbb{N}} \in l^2_r(\mathbb{N},\mathbb{C})$,
the function $v:=L(d)$ solves
$$\int_0^T v(t) e^{i \omega_k t} ds = d_k, \forall k \in \mathbb{N}.$$
\end{Cor}

\noindent \textbf{Proof of Corollary \ref{Cor:haraux1}:}
We define $\omega_{-k}:=-\omega_k, \forall k \in \mathbb{N}^*$.
Theorem \ref{thm:Haraux} ensures that the family
$(e^{i \omega_k t})_{k \in \mathbb{Z}}$ is a Riesz basis of
$F:=\text{Adh}_{L^2(0,T)} ( \text{Span} \{ e^{i \omega_k t} ; k \in \mathbb{Z} \})$.
Thanks to Proposition \ref{Riesz=solution-moment}, the map
$$\begin{array}{cccc}
J: & F & \rightarrow & l^2(\mathbb{Z},\mathbb{C}) \\
   & v & \mapsto     & 
\left( \int_0^T v(t) e^{i \omega_k t} dt \right)_{k \in \mathbb{Z}}
\end{array}$$
is an isomorphism. For $d=(d_k)_{k \in \mathbb{N}} \in l^2_r(\mathbb{N},\mathbb{C})$,
we define $\tilde{d}:=(\tilde{d}_k)_{k \in \mathbb{Z}} \in l^2(\mathbb{Z},\mathbb{C})$
by $\tilde{d}_k:=d_k$ if $k \geqslant 0$ and $\overline{d_{-k}}$ if $k<0$.
Now, we define 
$L:l^2_r(\mathbb{N},\mathbb{C}) \rightarrow L^2((0,T),\mathbb{R})$ by
$L(d)=J^{-1}(\tilde{d})$. The map $L$ takes values in real valued functions
because $\tilde{d}_{-k}=\overline{\tilde{d}_k}, \forall k \in \mathbb{N}$
for every $d \in l^2_r(\mathbb{N},\mathbb{C})$. $\Box$
\\

Theorem \ref{thm:Haraux} is also crucial in the proof of the
following statement, used in the proof of Proposition \ref{Control-Linearise}.

\begin{Cor} \label{Cor:haraux}
Let $T>0$ and $(\omega_{k})_{k \in \mathbb{N}}$ be an increasing sequence of $[0,+\infty)$
such that $\omega_{0}=0$ and 
\begin{eqnarray} \label{diffwkinfty}
\omega_{k+1} - \omega_{k} \rightarrow + \infty \text{ when } k \rightarrow + \infty.
\end{eqnarray}
There exists a continuous linear map
$$\begin{array}{cccc}
L: & \mathbb{R} \times l^2_r(\mathbb{N},\mathbb{C}) & \rightarrow & L^2((0,T),\mathbb{R}) \\
   &        (\tilde{d},d)                           & \mapsto     & L(\tilde{d},d)
\end{array}$$
such that, for every $\tilde{d} \in \mathbb{R}$,
$d=(d_k)_{k \in \mathbb{N}} \in l^2_r(\mathbb{N},\mathbb{C})$,
the function $v:=L(\tilde{d},d)$ solves
\begin{equation} \label{trigo-moment-pb}
\begin{array}{l}
\int_{0}^{T} v(t) e^{i \omega_{k} t} dt = d_{k} , \forall k \in \mathbb{N}, \\
\int_0^T t v(t) dt = \tilde{d}.
\end{array}
\end{equation}
\end{Cor}

\noindent \textbf{Proof of Corollary \ref{Cor:haraux}:} Let $\omega_k:=-\omega_{-k}$,
for every $k \in \mathbb{Z}$ with $k<0$. From Proposition \ref{thm:Haraux},
$\Theta:=(e^{i \omega_k t})_{k \in \mathbb{Z}}$ is a Riesz basis of
$\text{Adh}_{L^2(0,T)} ( \text{Span} \Theta)$.

\emph{First step: We prove that the family
$\widetilde{\Theta}:=\{t,e^{i \omega_k t}; k \in \mathbb{Z} \}$ is minimal in
$L^2(0,T)$.} 

Working by contradiction, we assume that 
$\widetilde{\Theta}$ is not minimal in $L^2(0,T)$. Then, necessarily
\begin{equation} \label{absurde-t}
t \in \text{Adh}_{L^2(0,T)} \text{Span} \Theta.
\end{equation}
With successive integrations, we get
$$t^j \in 
\text{Adh}_{C^{0}[0,T]} \Big( \text{Span}  \widetilde{\Theta} \Big),
\forall j \in \mathbb{N} \text{ with } j \geqslant 2.$$
The Stone Weierstrass theorem ensures that $\{1,t^j ; j \in \mathbb{N}, j \geqslant 2\}$
is dense in $C^0([0,T],\mathbb{C})$, thus, it is also dense in $L^2(0,T)$.
From (\ref{absurde-t}), we deduce that 
$\text{Span} \Theta$ is dense in
$L^2(0,T)$. This is a contradiction, because, 
thanks to Theorem \ref{thm:Haraux},
for every $\omega \in \mathbb{R}- \{ \omega_k, k \in \mathbb{Z} \}$,
the family $\{ e^{i \omega t}, e^{i \omega_k t}; k \in \mathbb{Z} \}$ 
is minimal, i.e.
$$e^{i \omega t} \notin \text{Adh}_{L^2(0,T)} \Big( \text{Span} \Theta \Big).$$

\emph{Second step: We conclude.}

For $k<0$, we define $d_k := \overline{d_{-k}}$.
Let $\{ \tilde{\xi}, \xi_k ; k \in \mathbb{Z} \}$ be the biorthogonal family
to $\{t, e^{i \omega_k t} ; k \in \mathbb{Z} \}$. From Theorem \ref{thm:Haraux},
there exists $C>0$ and a unique solution 
$v \in \text{Adh}_{L^2(0,T)} ( \text{Span} \Theta)$
of
$$\int_0^T v(t) e^{i \omega_k t} dt = d_k, \forall k \in \mathbb{Z}$$
and it satisfies
$$\|v\|_{L^2(0,T)} \leqslant C \left( \sum_{k \in \mathbb{Z}} |d_k|^2 \right)^{1/2}.$$
The uniqueness guarantees that $v$ is real valued. Let us define
$$L(\tilde{d},d):= u := v + \Big( \tilde{d}-\int_0^T t v(t) dt \Big) \tilde{\xi}.$$
Then, $u$ is real valued (because $v$ and $\tilde{\xi}$ are), 
$u$ solves (\ref{trigo-moment-pb}) and 
$$\begin{array}{ll}
\|u\|_{L^2} 
& 
\leqslant \| v \|_{L^2} +
\Big( |\tilde{d}| + \Big| \int_0^T t v(t) dt \Big| \Big) \| \tilde{\xi}\|_{L^2}
\\ & \leqslant
\| v \|_{L^2} \Big( 1 + \sqrt{\frac{T^3}{3}} \| \tilde{\xi}\|_{L^2} \Big)
+ |\tilde{d}| \| \tilde{\xi}\|_{L^2}
\\ & \leqslant
\left( C \Big( 1 + \sqrt{\frac{T^3}{3}} \| \tilde{\xi}\|_{L^2} \Big) 
+ \| \tilde{\xi}\|_{L^2} \right)
\left( |\tilde{d}|^2 + \sum_{k \in \mathbb{Z}} |d_k|^2 \right)^{1/2}. \Box
\end{array}$$

For the wave equation, the gap between two successive frequencies
does not tend to infinity, so we will need the following Corollary 
which is proved similarly.

\begin{Cor}\label{Corwave}
Let $T>2$. We make the same assumptions as Corollary \ref{Cor:haraux} except that we assume 
$$\omega_{k+1} - \omega_{k} \geqslant \pi$$
instead of (\ref{diffwkinfty}). Then, we have the same conclusion as Corollary \ref{Cor:haraux}. 
\end{Cor}

\begin{Cor} \label{Cor3}
Let $(\omega_{k})_{k \in \mathbb{N}}$ be an increasing sequence of $[0,+\infty)$
such that $\omega_{0}=0$ and 
$$\omega_{k+1} - \omega_{k} >\gamma >0.$$
There exists an nondecreasing function 
$$\begin{array}{cccc}
C: & [0,+\infty) & \rightarrow & \mathbb{R}^*_+ \\
   &     T       & \mapsto     & C(T)
\end{array}$$
such that, for every $T>0$ and for every $g \in L^2(0,T)$, we have
$$\left( \sum\limits_{k=0}^{\infty} 
\Big| \int_0^T g(t) e^{i \omega_k t} dt  \Big|^{2}
\right)^{1/2}
\leqslant C(T) \|g\|_{L^2(0,T)}.$$
\end{Cor}

\noindent \textbf{Proof of Corollary \ref{Cor3}:}
The existence of $C(T)$, for large $T\geqslant 2\pi/\gamma+1$,
is a consequence of Theorem \ref{thm:Haraux}  and 
Proposition \ref{Prop:Riesz} \textbf{(3)}. 
Let us choose for $C(T)$ the smallest value possible for this constant. For $T\leq 2\pi/\gamma+1 $, we choose $C(T)=C(2\pi/\gamma+1)$. 
Let $0<T_1<T_2<+\infty$,
$g \in L^2(0,T_1)$ and $\tilde{g} \in L^2(0,T_2)$ be defined by
$\tilde{g}=g$ on $(0,T_1)$ and $0$ on $(T_1,T_2)$. By applying the
inequality on $\tilde{g}$, we get $C(T_1) \leqslant C(T_2)$. $\Box$

\tableofcontents

\bibliography{biblio3CL}

\begin{thebibliography}{10}

\bibitem{Boscain-Adami}
R.~Adami and U.~Boscain.
\newblock {Controllability of the Schroedinger Equation via Intersection of
  Eigenvalues}.
\newblock {\em Proceedings of the 44rd IEEE Conference on Decision and Control
  December 12-15, 2005, Seville, (Spain). Also on 'Control Systems: Theory,
  Numerics and Applications, Roma, Italia 30 Mar - 1 Apr 2005, POS, Proceeding
  of science.}

\bibitem{Agrachev-Chambrion}
A.~Agrachev and T.~Chambrion.
\newblock An estimation of the controllability time for single-input systems on
  compact lie groups.
\newblock {\em ESAIM Control Optim. Calc. Var.}, 12(3):409--441, 2006.

\bibitem{Agrachev-book}
A.~Agrachev and Y.~L. Sachkov.
\newblock {\em Control theory from the geometric viewpoint.}
\newblock Encyclopaedia of Mathematical Sciences, 87. Control Theory and
  Optimization, II. Springer-Verlag, Berlin, 2004.

\bibitem{Agrachev9}
A.~Agrachev and A.~V. Sarychev.
\newblock {Navier-Stokes equations: controllability by means of low modes
  forcing}.
\newblock {\em J. Math. Fluid Mech.}, 7(1):108--152, 2005.

\bibitem{Albertini}
F.~Albertini and D.~D'Alessandro.
\newblock Notions of controllability for bilinear multilevel quantum systems.
\newblock {\em IEEE Transactions on Automatic Control}, 48(8):1399--1403, 2003.

\bibitem{Alinhac-Gerard}
S.~Alinhac and P.~Gérard.
\newblock {\em {Opérateurs pseudo-différentiels et théorème de Nash-Moser}}.
\newblock Intereditions (Paris), collection Savoirs actuels, 1991.

\bibitem{Altafini}
C.~Altafini.
\newblock Controllability of quantum mechanical systems by root space
  decomposition of su(n).
\newblock {\em J. Mathematical Physics}, 43(5):2051--2062, 2002.

\bibitem{AntonNLSbouleradial}
R.~Anton.
\newblock {Cubic nonlinear Schrödinger equation on three dimensional balls with
  radial data}.
\newblock {\em Comm. Partial Differential Equations}, 33(10-12):1862--1889,
  2008.

\bibitem{Avdonin}
S.A. Avdonin and S.A. Ivanov.
\newblock {\em Families of exponentials : the method of moments in
  controllability problems for distributed parameter systems}.
\newblock Cambridge New York , Cambridge University Press, 1995.

\bibitem{ball-marsden-slemrod}
J.M. Ball, J.E. Marsden, and M.~Slemrod.
\newblock Controllability for distributed bilinear systems.
\newblock {\em SIAM J. Control and Optim.}, 20, July 1982.

\bibitem{Baudouin}
L.~Baudouin.
\newblock {A bilinear optimal control problem applied to a time dependent
  Hartree-Fock equation coupled with classical nuclear dynamics}.
\newblock {\em Port. Math. (N.S.)}, 63(3):293--325, 2006.

\bibitem{Baudouin-Kavian-Puel}
L.~Baudouin, O.~Kavian, and J.-P. Puel.
\newblock {Regularity for a Schrödinger equation with singular potential and
  application to bilinear optimal control}.
\newblock {\em J. Diff. Eq.}, 216:188--222, 2005.

\bibitem{Baudouin-Salomon}
L.~Baudouin and J.~Salomon.
\newblock {Constructive solutions of a bilinear control problem for a
  Schrödinger equation}.
\newblock {\em Systems and Control Letters}, 57(6):453--464, 2008.

\bibitem{KB-JMPA}
K.~Beauchard.
\newblock {Local Controllability of a 1-D Schrödinger equation}.
\newblock {\em J. Math. Pures et Appl.}, 84:851--956, July 2005.

\bibitem{SchroLgVar}
K.~Beauchard.
\newblock {Controllability of a quantum particule in a 1D variable domain}.
\newblock {\em ESAIM:COCV}, 14(1):105--147, 2008.

\bibitem{KB-poutres}
K.~Beauchard.
\newblock {Local Controllability of a 1-D beam equation}.
\newblock {\em SIAM J. Control Optim.}, 47(3):1219--1273, 2008.

\bibitem{KB:SchroTIL}
K.~Beauchard.
\newblock {Local Controllability of a 1-D bilinear Schr\"odinger equation: a
  simpler proof}.
\newblock {\em (preprint)}, 2009.

\bibitem{KB-JMC}
K.~Beauchard and J.-M. Coron.
\newblock Controllability of a quantum particle in a moving potential well.
\newblock {\em J. Functional Analysis}, 232:328--389, 2006.

\bibitem{KB-JMC-MM-PR}
K.~Beauchard, J.-M. Coron, M.~Mirrahimi, and P.~Rouchon.
\newblock {Implicit Lyapunov control of finite dimensional Schrödinger
  equations}.
\newblock {\em System and Control Letters}, 56:388--395, 2007.

\bibitem{KB-MM}
K.~Beauchard and M.~Mirrahimi.
\newblock Practical stabilization of a quantum particle in a one-dimensional
  infinite square potential well.
\newblock {\em SIAM J. Contr. Optim.}, 48(2):1179--1205, 2009.

\bibitem{Boscain-Chambrion-Charlot}
U.~Boscain, T.~Chambrion, and G.~Charlot.
\newblock Nonisotropic 3-level quantum systems: complete solutions for minimum
  time and minimum energy.
\newblock {\em Discrete Contin. Dyn. Syst. Ser. B}, 5(4):957--990, 2005.

\bibitem{Boscain-Charlot-2005}
U.~Boscain and G.~Charlot.
\newblock Resonance of minimizers for $n$-level quantum systems with an
  arbitrary cost.
\newblock {\em ESAIM Control Optim. Calc. Var.}, 10(4):593--614, 2004.

\bibitem{Boscain-Charlot-Gauthier-book}
U.~Boscain, G.~Charlot, and J.-P. Gauthier.
\newblock {\em {Optimal control of the Schr\"odinger equation with two or three
  levels, Nonlinear and adaptive control (Sheffield, 2001)}}, volume 281.
\newblock Lecture Notes in Control and Inform. Sci., Springer, Berlin, 2003.

\bibitem{Boscain-Charlot-Gauthier}
U.~Boscain, G.~Charlot, J.-P. Gauthier, S.~Gu\'erin, and H.-R. Jauslin.
\newblock Optimal control in laser-induced population transfer for two- and
  three-level quantum systems.
\newblock {\em J. Math. Phys.}, 43(5), 2002.

\bibitem{Gauthier-Jauslin}
U.~Boscain, G.~Charlot, J.-P. Gauthier, S.~Gu\'erin, and H.-R. Jauslin.
\newblock Optimal control in laser-induced population transfer for two- and
  three-level quantum systems.
\newblock {\em J. Math. Phys.}, 43:2107--2132, 2002.

\bibitem{Boscain-Mason}
U.~Boscain and P.~Mason.
\newblock Time minimal trajectories for a spin 1/2 particle in a magnetic
  field.
\newblock {\em J. Math. Phys.}, 47(6):062101--29, 2006.

\bibitem{Brockett}
R.~Brockett.
\newblock Lie theory and control systems defined on spheres.
\newblock {\em SIAM J. Appl. Math.}, 25(2):213--225, 1973.

\bibitem{Burq}
N.~Burq.
\newblock Contrôle de l'équation des plaques en présence d'obstacles
  strictement convexes.
\newblock {\em M\'emoire de la S.M.F.}, 55, 1993.

\bibitem{Crepeau_Cerpa}
E.~Cerpa and E.~Crépeau.
\newblock Boundary controlability for the non linear korteweg-de vries equation
  on any critical domain.
\newblock {\em Ann. IHP, Analyse Non Linéaire (in press)}.

\bibitem{Chambrion-et-al}
T.~Chambrion, P.~Mason, M.~Sigalotti, and M.~Boscain.
\newblock {Controllability of the discrete-spectrum Schrödinger equation driven
  by an external field.}
\newblock {\em Ann. Inst. H. Poincaré Anal. Non Linéaire}, 26(1):329--349,
  2009.

\bibitem{jmc-tank}
J.-M. Coron.
\newblock {Local Controllability of a 1-D Tank Containing a Fluid Modeled by
  the shallow water equations}.
\newblock {\em ESAIM: COCV}, 8:513--554, June 2002.

\bibitem{JMC-CRAS-Tmin}
J.-M. Coron.
\newblock On the small-time local controllability of a quantum particule in a
  moving one-dimensional infinite square potential well.
\newblock {\em C. R. Acad. Sciences Paris, Ser. I}, 342:103--108, 2006.

\bibitem{JMC-book}
J.-M. Coron.
\newblock {\em Control and nonlinearity}, volume 136.
\newblock Mathematical Surveys and Monographs, 2007.

\bibitem{Dehman-Gerard}
B.~Dehman, P.~Gérard, and G.~Lebeau.
\newblock {Stabilization and control for the nonlinear Schrödinger equation on
  a compact surface}.
\newblock {\em Mathematische Zeitschrift}, 254(4):729--749, December 2006.

\bibitem{CLB-Cances-Pilot}
{E. Cancès and C. Le Bris and M. Pilot}.
\newblock {Contrôle optimal bilinéaire d'une équation de Schrödinger}.
\newblock {\em CRAS Paris}, {330}:567--571, 2000.

\bibitem{Ervedoza-Puel}
S.~Ervedoza and J.-P. Puel.
\newblock Approximate controllability for a system of schrödinger equations
  modeling a single trapped ion,.
\newblock {\em Annales de l'Institut Henri Poincaré : Analyse non linéaire (to
  appear)}, 2009.

\bibitem{Fabre}
C.~Fabre.
\newblock {Résultats de contrôlabilité exacte interne pour l'équation de
  Schrödinger et leurs limites asymptotiques: application à certaines équations
  de plaques vibrantes. (French) [Results on exact internal controllability for
  the Schrödinger equation and their asymptotic limits: application to some
  vibrating-plate equations]}.
\newblock {\em Asymptotic Analalysis}, 5(4):343--379, 1992.

\bibitem{Haraux}
A.~Haraux.
\newblock {Séries lacunaires et contrôle semi-interne des vibrations d'une
  plaque rectangulaire}.
\newblock {\em J. Math. Pures et Appl.}, 68:457--465, 1989.

\bibitem{horm}
{H\"ormander}.
\newblock {On the Nash-Moser Implicit Function Theorem}.
\newblock {\em Annales Academiae Scientiarum Fennicae}, pages 255--259, 1985.

\bibitem{Teismann-et-al}
R.~Ilner, H.~Lange, and H.~Teismann.
\newblock {Limitations on the control of Schrödinger equations}.
\newblock {\em ESAIM:COCV}, 12(4):615--635, 2006.

\bibitem{Khaneja-Glaser-Brockett}
N.~Khaneja, S.~J. Glaser, and R.~Brockett.
\newblock Sub-riemannian geometry and time optimal control of three spin
  systems: quantum gates and coherence transfer.
\newblock {\em Phys. Rev. A (3)}, 65:032301, 11, 2002.

\bibitem{Teismann-et-alBIS}
H.~Lange and H.~Teismann.
\newblock {Controllability of the nonlinear Schrödinger equation in the
  vicinity of the ground state}.
\newblock {\em Math. Methods Appl. Sci.}, 30(13):1483--1505, 2007.

\bibitem{Lasiecka-Triggiani}
I.~Lasiecka and R.~Triggiani.
\newblock {Optimal regularity, exact controllability and uniform stabilization
  of Schrödinger equations with Dirichlet controls}.
\newblock {\em Differential and Integral Equations}, 5:571--535, 1992.

\bibitem{Lasiecka-Triggiani-Zhang}
I.~Lasiecka, R.~Triggiani, and X.~Zhang.
\newblock {Global uniqueness, observability and stabilization of
  nonconservative Schrödinger equations via pointwise Carlemann estimates}.
\newblock {\em J. Inverse Ill Posed-Probl.}, 12:183--231, 2004.

\bibitem{Laurentdim3}
C.~Laurent.
\newblock {Global controllability and stabilization for the nonlinear
  Schr\"odinger equation on some compact manifolf of dimension 3}.
\newblock {\em submitted}, 2009.

\bibitem{Laurent}
C.~Laurent.
\newblock {Global controllability and stabilization for the nonlinear
  Schrödinger equation on an interval}.
\newblock {\em ESAIM:COCV (to appear)}, 2009.

\bibitem{lebeau}
G.~Lebeau.
\newblock {Contrôle de l'équation de Schrödinger}.
\newblock {\em J. Math. Pures Appl.}, 71:267--291, 1992.

\bibitem{Machtyngier}
Machtyngier.
\newblock {Exact controllability for the Schrödinger equation}.
\newblock {\em SIAM J. Contr. Opt.}, 32:24--34, 1994.

\bibitem{MM}
M.~Mirrahimi.
\newblock Lyapunov control of a quantum particle in a decaying potential.
\newblock {\em {Ann. Inst. H. Poincaré (c) Nonlinear Analysis}}, 26:1743--1765,
  2009.

\bibitem{HarmOsc}
M.~Mirrahimi and P.~Rouchon.
\newblock Controllability of quantum harmonic oscillators.
\newblock {\em IEEE Trans. Automatic Control}, 49(5):745--747, 2004.

\bibitem{MM-PR-GT}
M.~Mirrahimi, P.~Rouchon, and G.~Turinici.
\newblock {Lyapounov control of bilinear Schrödinger equations}.
\newblock {\em Automatica}, 41:1987--1994, 2005.

\bibitem{Nersesyan2}
V.~Nersesyan.
\newblock {Global approximate controllability for Schrödinger equation in
  higher Sobolev norms and applications}.
\newblock {\em (preprint)}, 2009.

\bibitem{Nersesyan1}
V.~Nersesyan.
\newblock {Growth of Sobolev norms and controllability of Schrödinger
  equation}.
\newblock {\em Comm. Math. Phys. (to appear)}, 2009.

\bibitem{ramakrishna-et-al-95}
V.~Ramakrishna, M.~Salapaka, M.~Dahleh, and H.~Rabitz.
\newblock {Controllability of molecular systems}.
\newblock {\em Phys. Rev. A}, 51(2):960--966, 1995.

\bibitem{Tucsnak}
K.~Ramdani, T.~Takahashi, G.~Tenenbaum, and M.~Tucsnak.
\newblock A spectral approach for the exact observability of
  infinite-dimensional systems with skew-adjoint generator.
\newblock {\em Journal of Functional Analysis}, 226:193--229, 2005.

\bibitem{RosierKdV}
L.~Rosier.
\newblock Exact boundary controllability for the korteweg-de vries equation on
  a bounded domain.
\newblock {\em ESAIM: Control, Optimisation and Calculus of Variations},
  (2):33--55, 1997.

\bibitem{Rosier-Zhang}
L.~Rosier and B.-Y. Zhang.
\newblock {Local exact controllability and stabilizability of the nonlinear
  Schrödinger equation on a bounded interval.}
\newblock {\em SIAM J. Control Optim.}, 48(2):972--992, 2009.

\bibitem{rouchon-oct2002}
P.~Rouchon.
\newblock Control of a quantum particule in a moving potential well.
\newblock {\em 2nd IFAC Workshop on Lagrangian and Hamiltonian Methods for
  Nonlinear Control, Seville}, 2003.

\bibitem{Russel_Zhang428}
D.~Russel and B.-Y. Zhang.
\newblock {Exact controllability and stabilizability of the Korteweg-de Vries
  equation.}
\newblock {\em Trans. Amer. Math. Soc.}, 348(9):3643--3672, 1996.

\bibitem{Shirikyan446}
A.~Shirikyan.
\newblock {Approximate controllability of three-dimensional Navier-Stokes
  equations}.
\newblock {\em Comm. Math. Phys.}, 266(1):123--151, 2006.

\bibitem{Sussmann-Jurdjevic}
H.J. Sussmann and V.~Jurdjevic.
\newblock Controllability of nonlinear systems.
\newblock {\em J. Differential Equations}, 12:95--116, 1972.

\bibitem{Defranceschi-LeBris}
G.~Turinici.
\newblock On the controllability of bilinear quantum systems.
\newblock {\em {In C. Le Bris and M. Defranceschi, editors, Mathematical Models
  and Methods for Ab Initio Quantum Chemistry}}, {volume 74 of Lecture Notes in
  Chemistry, Springer}, 2000.

\bibitem{turinici-rabitz-01}
G.~Turinici and H.~Rabitz.
\newblock Quantum wave function controllability.
\newblock {\em Chem. Phys.}, 267:1--9, 2001.

\bibitem{Zhang_513}
B.-Y. Zhang.
\newblock {Exact boundary controllability of the Korteweg-de Vries equation}.
\newblock {\em SIAM J. Cont. Optim.}, 37(2):543--565, 1999.

\bibitem{zuazua-wave}
E.~Zuazua.
\newblock Exact controllability for semilinear wave equations in one space
  dimension.
\newblock {\em Ann. Inst. H. Poincaré Anal. Non Linéaire}, 10(1):109--129,
  1993.

\end{thebibliography}
\bibliographystyle{plain}

\end{document}